\documentclass[reqno,12pt]{amsart}
\usepackage{enumitem,stix}
\usepackage{mathtools}
\usepackage{epsf}
\usepackage{mathtools}
\usepackage{ytableau}
\usepackage{todonotes}
\usepackage{verbatim}
\usepackage{amsmath}
\usepackage{amsopn}
\usepackage{amssymb}
\usepackage{latexsym}
\usepackage{amsfonts}
\usepackage{amsthm}
\usepackage{tikz}
\usepackage[foot]{amsaddr}

\newtheorem{thm}{Theorem}[section]
\newtheorem{prop}[thm]{Proposition}

\newtheorem{defn}[thm]{Definition}

\newtheorem{obs}[thm]{Observation}

\newtheorem{lem}[thm]{Lemma}
\newtheorem{exa}[thm]{Example}

\title{Geometric view of interval posets of permutations}

\author{Eli Bagno, Estrella Eisenberg, Shulamit Reches, and Moriah Sigron}\address{\textup{Eli Bagno, Shulamit Reches and Moriah Sigron  are members of the Department of Mathematics in Jerusalem College of Technology and
the Department of Mathematics in the  Michlala College of Jerusalem.}}
\address{\textup{Estrella Eisenberg is a member of the Department of Mathematics, Jerusalem College of Technology.}}

\begin{document}
\begin{abstract}
In a recent study, Tenner introduced the concept of the interval poset of a permutation to effectively represent all intervals and their inclusions within a permutation. In this paper, we present a new geometric viewpoint on interval posets. We establish a one-to-one correspondence between the set of interval posets for permutations of size $n$ and a specific subset of dissections of a convex polygon with $n+1$ sides. Through this correspondence, we investigate various intriguing subsets of interval posets and uncover their connections with specific polygon dissections.
\end{abstract}
\maketitle
\section{Introduction}
    The notion of an interval poset for a permutation, initially introduced by Tenner in \cite{T}, effectively represents all intervals within a permutation and their inclusion relationships comprehensively.
Before providing the formal definition, we recall some preliminary notions.

\begin{defn}
Let $\mathcal{S}_n$ be the symmetric group on $n$ elements and let $\pi\in \mathcal{S}_n$.
An interval in $\pi$ is a set of consecutive values that appear in consecutive positions of $\pi$.
Clearly, $[n]:=[1,n]$ is an interval, as well as $\{i\}$ for each $i \in [n]$. These are called {\it trivial intervals}. The other intervals are called {\it proper}.
\end{defn}
For example, for the permutation $\pi=314297856$ all the proper intervals are: $[5,9]=\{5,6,7,8,9\}$, $[1,4],[5,6],[7,8],[7,9]$ and $[5,8]$.
A permutation $\pi \in \mathcal{S}_n$ is called {\em simple} if it does not have any proper interval; for example, the permutation $3517246$ is simple.
Following Tenner \cite{T}, we define the {\it interval poset of a permutation} as follows.
\begin{figure}[!ht]
     \centering
           \includegraphics[scale=0.30]{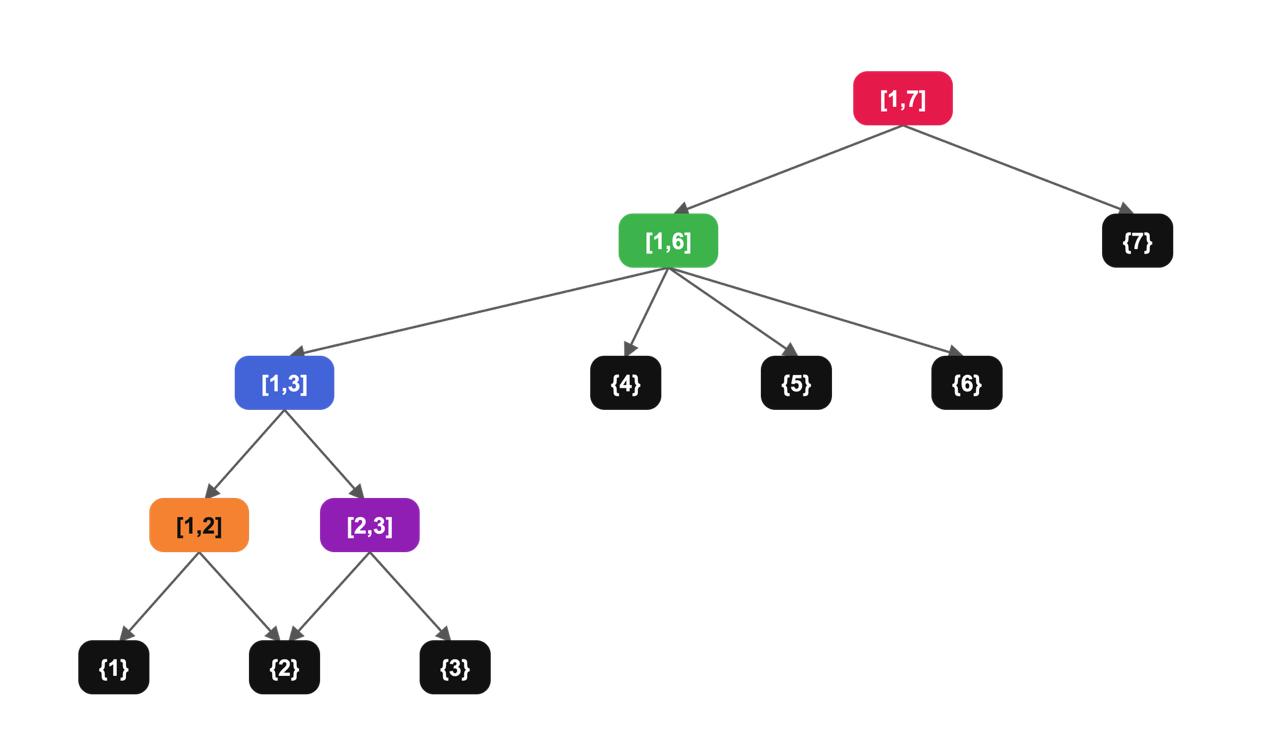}
\caption{Interval poset of the permutations: 5123647, 5321647, 4612357, 4632157, 7463215, 7461235, 7532164, 7512364.}
 \label{shulamit}
 \end{figure}
\begin{figure}[!ht]
     \centering
    \includegraphics[scale=0.20]{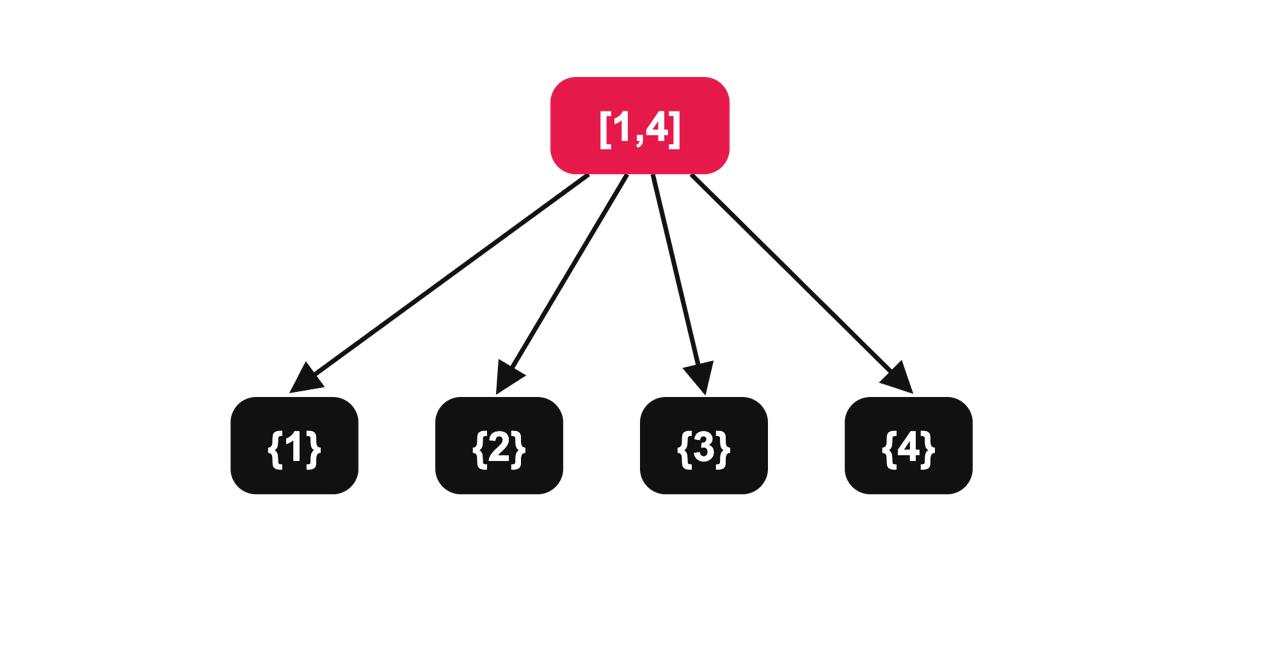}
\caption{Interval poset of permutations 3142 and 2413.  }
 \label{eli}
 \end{figure}

\begin{defn} \label{def interval poset}
The {\it interval poset of a permutation} $\pi \in \mathcal{S}_n$ is the poset $P(\pi)$ whose elements are the non-empty intervals of $\pi$; the order is defined by set inclusion.
The minimal elements are the intervals of size $1$. See, for example, Fig. \ref{shulamit} and Fig. \ref{eli}.
\end{defn}
In \cite{T}, the interval poset is arranged in the plane so that the elements covered by each node are positioned in increasing order from left to right, though \cite{BCI} offers an alternative embedding of the same poset.

 An interval poset might correspond to multiple permutations; note that all simple permutations of a given size $n$ have the same interval poset. See an example in Fig. \ref{eli}.
The aspects of interval posets of permutations that have already been investigated are structural and enumerative, see \cite{T} and \cite{BCI}.
Bouvel, Cioni and Izart \cite{BCI} provided a formula for the number of interval posets with $n$ minimal elements, i.e., arising from permutations of $n$ elements (sequence A348479 in \cite{OEIS}).
Tenner enumerated binary interval posets, i.e., posets in which each vertex covers at most two elements.
These are exactly the interval posets of the separable permutations (i.e., avoiding both  $2413$ and $3142$). These posets are counted by the small Schröder numbers, see Corollary 6.3 in \cite{T}.
Tenner also showed that the number of binary tree interval posets is $2C_{n-1}$, where $C_n$ denotes the $n-$th Catalan number. Here, a tree interval poset is an interval poset whose Hasse diagram is a tree.

Subsequently, Bouvel, Cioni, and Izart in \cite{BCI} enumerated the tree interval posets (sequence A054515 in \cite{OEIS}).
 This sequence also enumerates the number of ways to place
non-crossing diagonals in a convex $(n+1)$-gon such that no quadrilateral is
created. This coincidence raises the question of whether
there exists a direct bijection between these two classes.
The present work originated from an investigation of this question.
Rather than constructing a bijection only for tree interval posets, we show
that the natural geometric setting is considerably richer. By introducing a
natural class of polygon dissections, called \emph{diagonally framed
dissections}, we obtain a bijection with the entire class of interval posets.
Our main result, stated in Theorem~\ref{bijection}, establishes this
correspondence. Besides providing a geometric interpretation of interval
posets, it serves as the basis for new bijective proofs of several previously
known results, as well as for new enumerative results.

The existing enumerations of interval posets are obtained by
algebraic methods, primarily through generating functions; in contrast, our
approach is entirely bijective and geometric.
Restricting the bijection to tree interval posets immediately yields a
bijective proof of the enumeration obtained by Bouvel, Cioni and
Izart~\cite{BCI}. We further show that a suitable restriction of the same
bijection provides a geometric interpretation for interval posets arising
from block-wise simple permutations, introduced in~\cite{BERS}. Finally, for
separable permutations, we construct a different bijection between their
interval posets and non-crossing polygon dissections. This construction builds
on the classical bijection between binary trees and polygon triangulations.
\section{Preliminaries on Interval posets}
In this section, we gather some simple facts and observations on interval posets that will be used in the sequel. Most of the material presented in this section has already appeared elsewhere and is included here for the sake of self-containment.

Recall that in a poset $(P,\leq)$, $a$ {\it covers} $b$ if $b\leq a$ and there is no $c$ other than $a$ and $b$ such that $b \leq c \leq a$.

\begin{figure}
    \centering
    \includegraphics[scale=0.45]{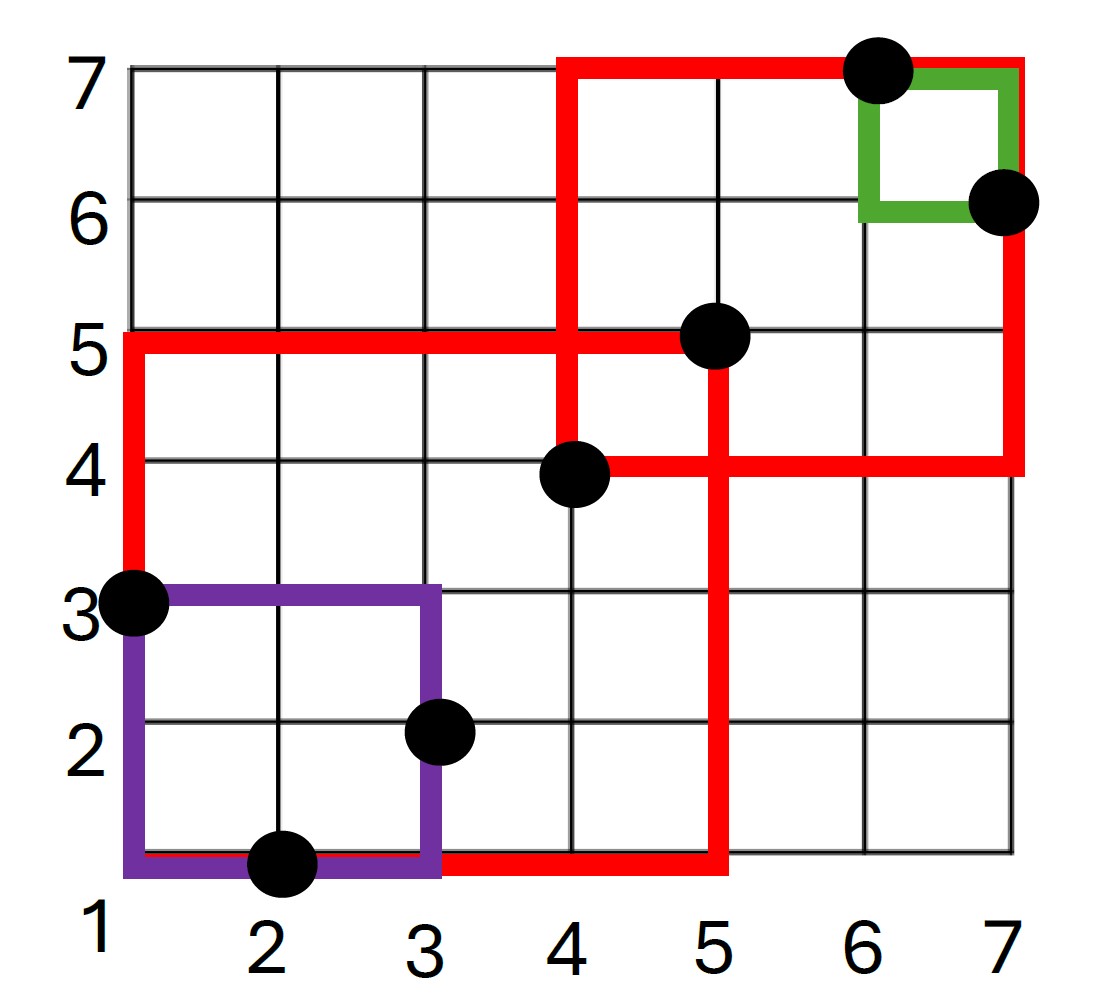}
    \caption{The permutation $\pi=3124576$ and some of its blocks in a “standard permutation plot.”}
    \label{f3124576}
\end{figure}

A version of the following observation has appeared in \cite{T} as Lemma 3.1:
\begin{obs}\label{structure of interval posets 1}
Let $\pi \in \mathcal{S}_n$. If $I$ and $J$ are intervals of $\pi$ such that $I$ is not contained in $J$  and $J$ is not contained in $I$, and $I\cap J\neq \emptyset$, then $I\cap J$, $I\cup J$, $I \setminus J$ and $J\setminus I$ are intervals of $\pi$.
\end{obs}
For example, let $\pi=3124576$. Then $I=[1,5]$ and $J=[4,7]$ are intersecting intervals of $\pi$ and thus $I\cup J=[1,7],I\cap J=[4,5], I \setminus J=[1,3], J \setminus I=[6,7]$ are also intervals of $\pi$, as illustrated in Fig. \ref{f3124576}, which depicts the permutation $\pi$ in the common graphical way, see \cite{BE}.

The following observation appears as Corollary 3.8 in \cite{T}.
\begin{prop}\label{structure of interval posets 2}
If $P(\pi)$ is the interval poset of $\pi\in \mathcal{S}_n$, then no element of $P(\pi)$ has exactly $3$ elements covered by it.
\end{prop}
\begin{proof}
Suppose that an interval $I$ covers exactly three intervals
$I_1,I_2,I_3$.
First assume that two of these intervals intersect, say $I_1\cap I_2\neq\varnothing$.
Since neither interval contains the other (otherwise, one of them could not be covered by
$I$), Observation \ref{structure of interval posets 1} implies that $I_1\cup I_2$ is also an interval.
Moreover,
\[
I_1,I_2\subsetneq I_1\cup I_2\subsetneq I,
\]
contradicting the fact that $I$ covers $I_1$ and $I_2$.
Hence, the intervals $I_1,I_2,I_3$ are pairwise disjoint. Since they are
maximal proper subintervals of $I$, they partition $I$ into three consecutive
blocks. Contracting each block to a single point yields a permutation of size
$3$. Every permutation of size $3$ contains an interval of size $2$. Lifting
this interval back to $I$ gives a proper interval strictly contained in $I$
and strictly containing two of the intervals $I_1,I_2,I_3$, again
contradicting the assumption that $I$ covers them.
Therefore, no interval in $P(\pi)$ can cover exactly three elements.
\end{proof}
Tenner \cite{T} gave a characterization of interval posets, which we cite here right after the following definition.

\begin{figure}
    \centering
    \includegraphics[scale=0.45]{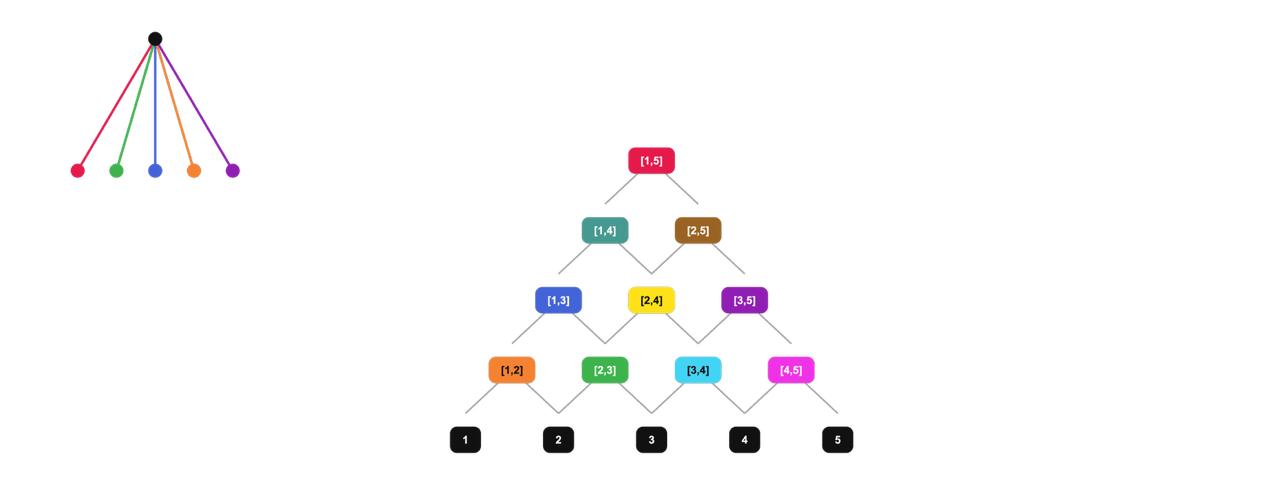}
    \caption{From left to right: a dual claw poset, an argyle poset.}
    \label{Argyle}
\end{figure}

\begin{defn}
\label{dual claw argyle}
\begin{enumerate}
\item A {\em dual claw poset} is a poset that has a unique maximal element and $k \geq 4$ minimal elements covered by that maximal element.
\item  An {\em argyle poset} of order $n$ is a poset which is isomorphic to the poset of all integer intervals of the set $\{1,\dots,n\}$, including the singletons $\{1\},\dots,\{n\}$,  with respect to set inclusion.
\end{enumerate}
 Clearly, the size of an argyle poset of order $n$ is $\binom{n+1}{2}$.
See Fig. \ref{Argyle}.
\end{defn}
Note that the interval poset of the identity permutation $1\cdots n$ is an argyle poset of order $n$.
Moreover, the interval poset of every simple permutation of order $n$ is the dual claw poset of order $n$.
We now quote Theorem 4.8 from \cite{T}, which explains how to construct a general interval poset, using the elementary building blocks defined above.
\begin{thm}
\label{Interval}
Interval posets are exactly those that can be constructed by starting with the
1-element poset, and recursively replacing minimal elements with dual claw posets and argyle posets.
\end{thm}
\section{Interval posets and convex polygons}
\label{the bijection}
\begin{figure}
    \centering
   \includegraphics[scale=0.35]{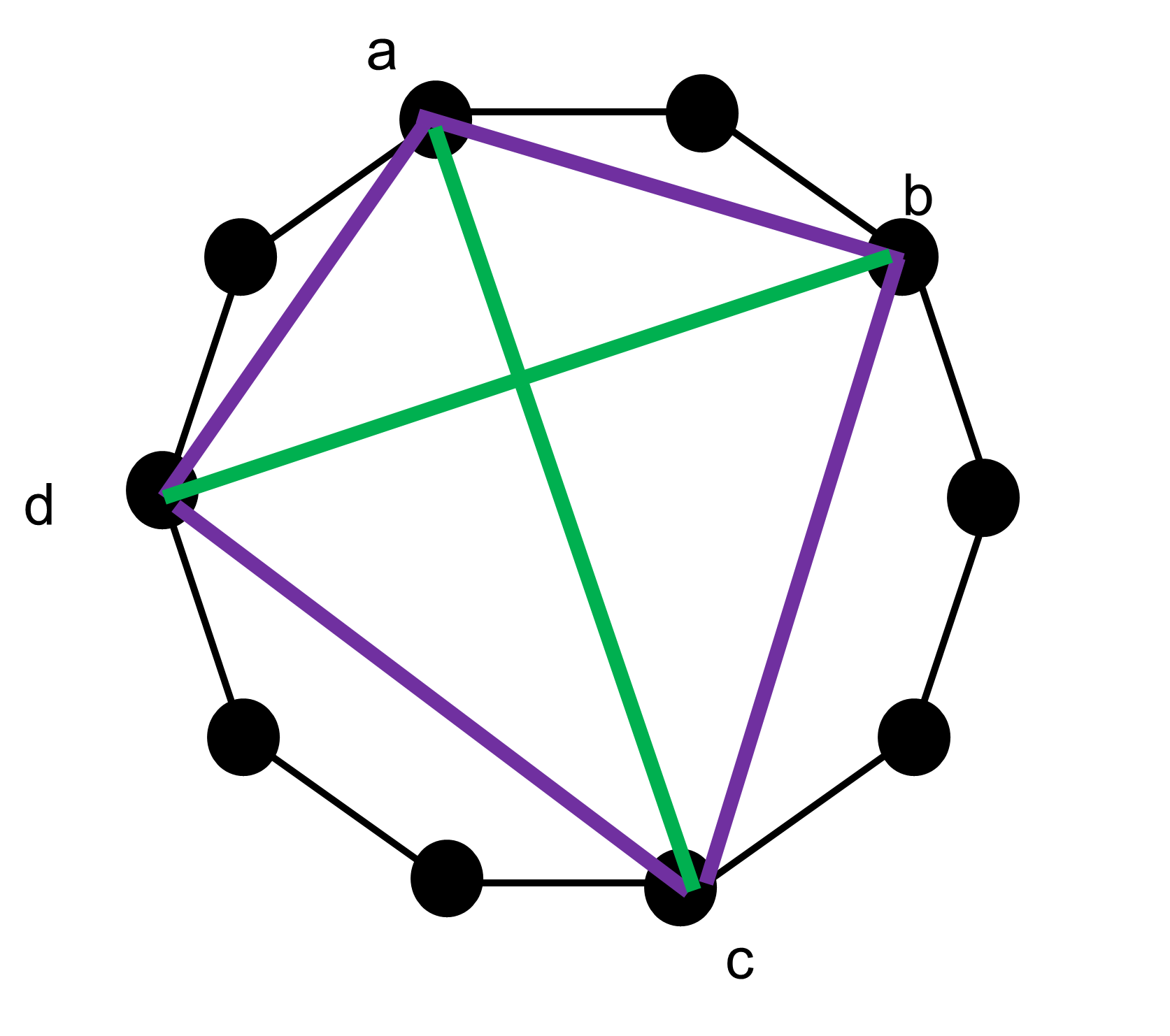}
    \caption{A $10$-gon with crossing diagonals $(a,c)$ and $(b,d)$ and their \emph{frame}.}
    \label{frame1}
\end{figure}

The main tool of this work is a bijection between the set of interval posets of permutations of size $n$ and a certain subset of the set of dissections of the convex $(n+1)$-gon, which we define below.  We show how to use this bijection as well as some restrictions of it to enumerate various subsets of the set of interval posets.
We identify a convex polygon with its set of vertices and denote a diagonal or an outer edge of the polygon
from vertex $i$ to vertex $j$ by $(i,j)$.
We start with the following well-known definition:
\begin{defn}
Let $P$ be a convex $n$-gon. A {\em dissection} of $P$ is a way of placing diagonals in the polygon.

Although we distinguish between the original (outer) edges of the polygon and its diagonals, occasionally, we also refer to diagonals as edges.

We use the term {\em crossing} to refer to edges that intersect at an interior point.

We use the term {\it quadrilateral} to mean a subpolygon with four vertices, whose four boundary edges belong to the dissection, but neither of the two diagonals joining opposite vertices belongs to the dissection (see the quadrilaterals $\{1,2,3,4\}$, $\{1,4,6,8\}$,$\{1,8,9,10\}$, $\{1,10,11,12\}$ in Fig.~\ref{example for bijection separable}).

Finally, for a given dissection of the polygon, we define the {\em graph of the dissection} to be the graph whose vertex set is the set of vertices of the polygon and whose edge set is the union of the outer edges and the diagonals of the dissection.
\end{defn}
We introduce the following definition, which will be used extensively throughout this paper.
\begin{defn}\label{def df}
A dissection of an $n$-gon will be called {\em diagonally framed} if for each pair of crossing diagonals, all their vertices are connected. Explicitly, if $(a,c)$ and $(b,d)$ are two crossing diagonals,
then the diagonals or outer edges $(a,b),(b,c),(c,d),(d,a)$ must also be part of the dissection (see Fig. \ref{frame1} for an example).
 \end{defn}
We define a function  $\Phi$ from the set of interval posets with $n$ minimal elements to the set of diagonally framed dissections of a convex $(n+1)$-gon without quadrilaterals and prove that it is a bijection.

\begin{figure}[!ht]
    \centering
\includegraphics[scale=0.20]{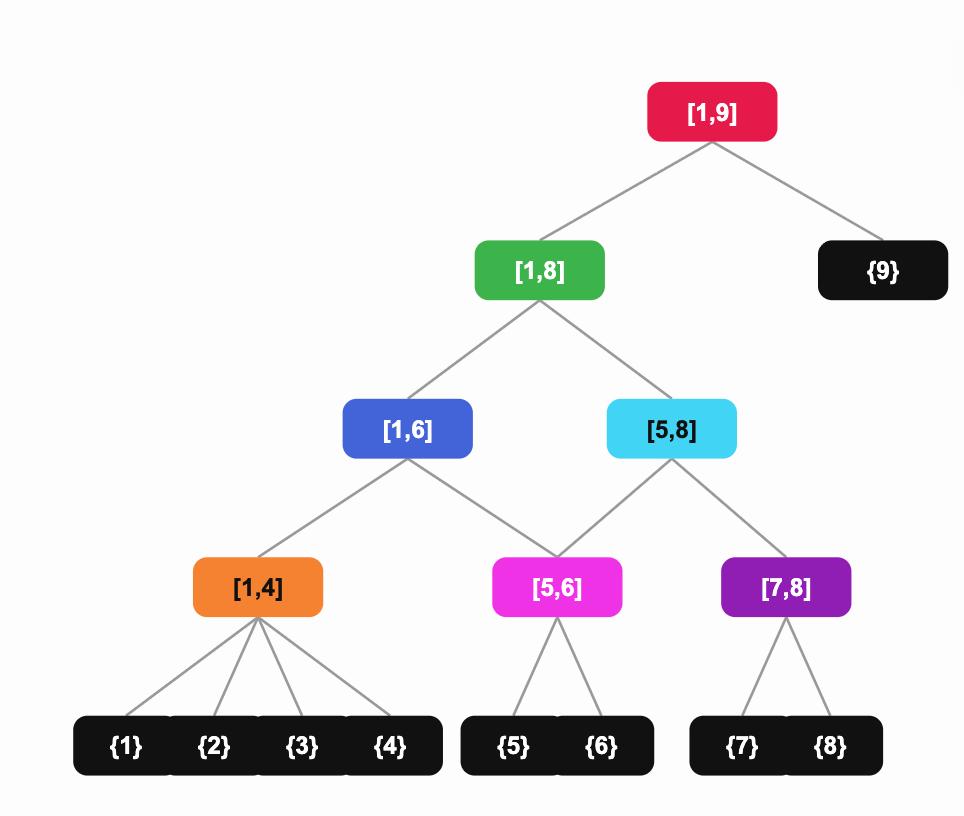}
\includegraphics[scale=0.30]{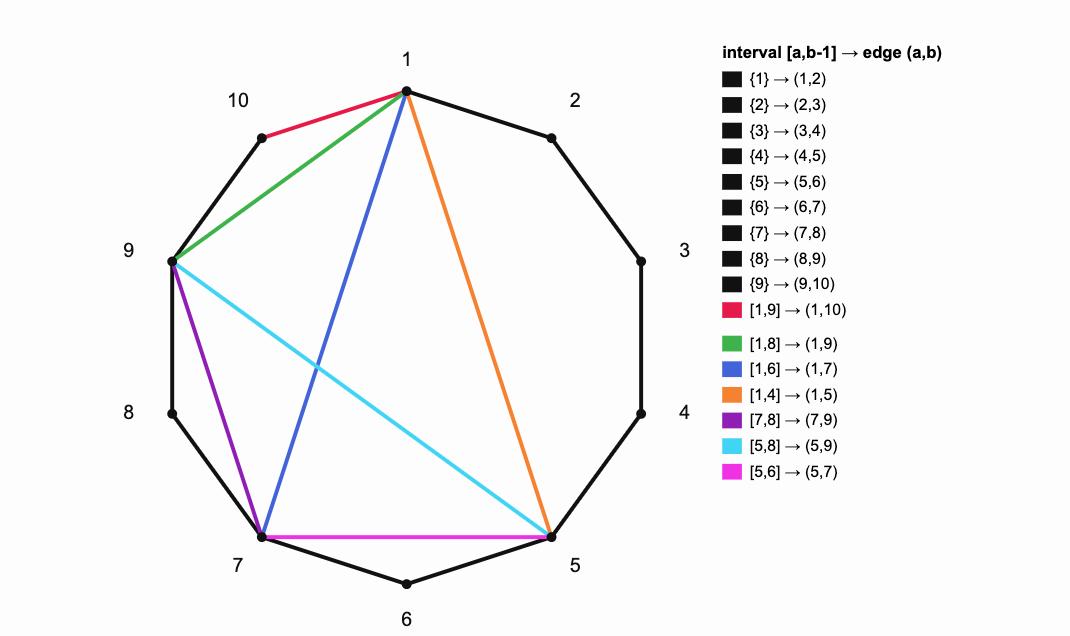}
    \caption{Up: the interval poset $P$. Down: the polygon $\Phi(P)$.}
    \label{polygon to poset}
\end{figure}
Let $P$ be the interval poset of some $\pi\in \mathcal{S}_n$. We set $\Phi(P)$ to be the dissection of the convex $(n+1)$-gon whose set of diagonals is $$\{ (a,b+1)|[a,b] \text{ is a non-minimal element of } P\},$$ i.e., to each non-singleton interval of the form $[a,b]$ corresponds a diagonal $(a,b+1)$ in $\Phi(P)$; note that singleton intervals of the form $
\{a\}$ correspond to outer edges $(a,a+1)$ in the polygon  (see  Fig. \ref{polygon to poset}).

\begin{figure}
   \centering
  \includegraphics[scale=0.15]{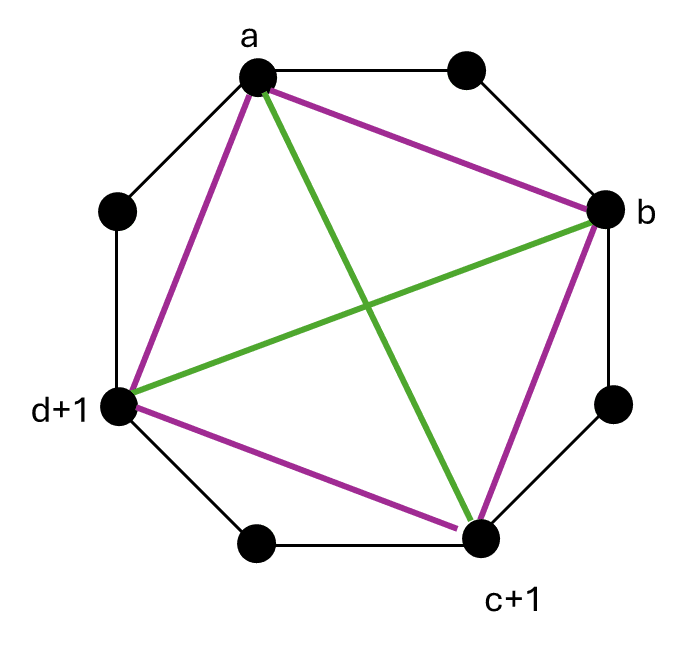}
    \caption{An octagon with crossing diagonals.} 
   \label{frame2}
\end{figure}

\begin{lem}\label{Phi(P) is diagonally framed}
    Let $P$ be an interval poset with $n$ minimal elements. Then the $(n+1)$-gon $\Phi(P)$ is diagonally framed.
\end{lem}
\begin{proof}
 If $(a,c+1)$ and $(b,d+1)$ are two crossing diagonals in $\Phi(P)$, where $a<b\leq c<d$, then $I=[a,c]$ and $J=[b,d]$ are intersecting intervals in $P$ and by Observation \ref{structure of interval posets 1}, we have that $I\cup J=[a,d]$, $I\cap J=[b,c]$, $I \setminus J=[a,b-1]$ and $J \setminus I=[c+1,d]$ are intervals in $P$ corresponding respectively to the diagonals or outer edges $(a,d+1), (b,c+1),(a,b)$ and $(c+1,d+1)$ (see Fig.  \ref{frame2}). \end{proof}

\begin{obs}
\label{argyle}
The interval poset $P$ is an argyle poset if and only if the graph of the dissection $\Phi(P)$ is complete.  Moreover, the interval poset $P$ is a dual claw poset if and only if the graph of $\Phi(P)$ has no internal edge, i.e., $\Phi(P)$ is the empty dissection.
\end{obs}

The following definitions and lemmas prepare the proof of Theorem \ref{bijection}. They show that every diagonally framed dissection is composed of complete subpolygons and empty subpolygons of size greater than four. 

\begin{figure}
    \centering
   \includegraphics[scale=0.39]{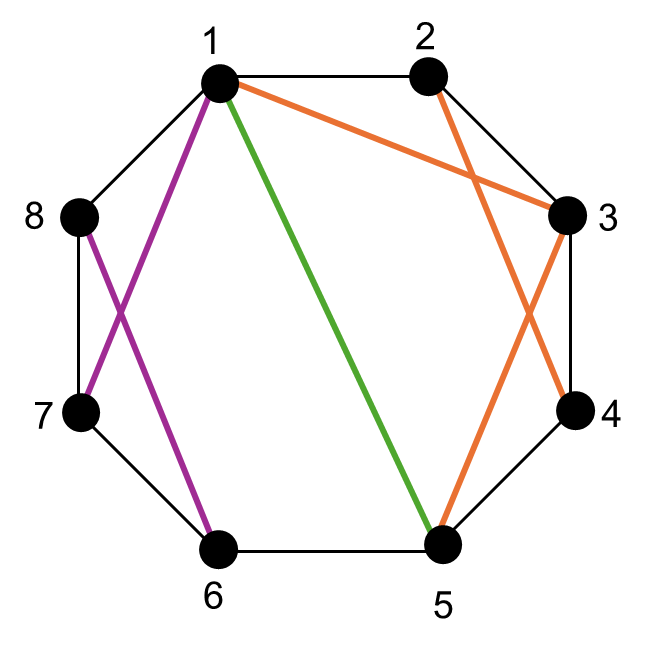}
   \caption{Dissection of $C_8$.}
  \label{Intersectional components}
\end{figure}

\begin{figure}
    \centering
   \includegraphics[scale=0.19]{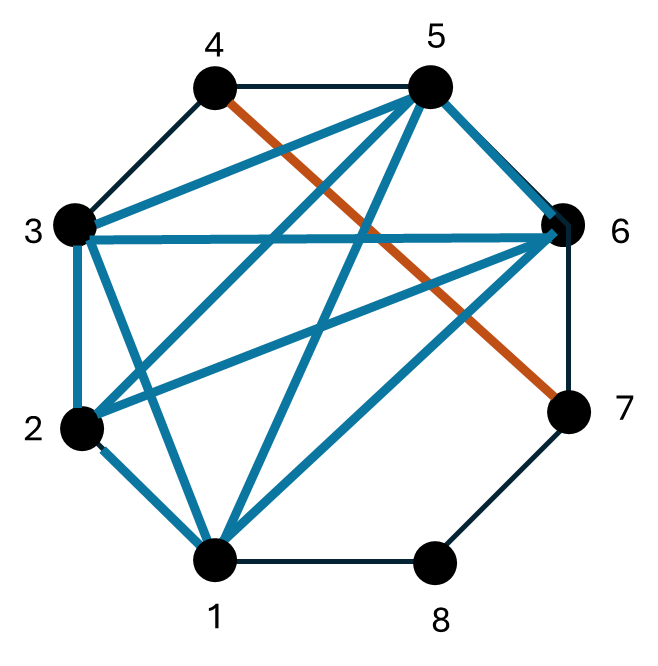}
   \caption{A Pentagon with a crossing diagonal.}
  \label{hexagon}
\end{figure}

\begin{defn}
For any dissection of a convex $n$-gon with set of diagonals $D=\{d_1,\dots,d_k\}$, define an equivalence relation on $D$ by $d_i \sim d_j$ if there is a sequence $(d_i=d_{i_1},\dots,d_{i_\ell}=d_j)$ such that for each $1\leq u \leq \ell-1$, $d_{i_u}$ crosses $d_{i_{u+1}}$ (i.e., intersects it in its interior).
Each equivalence class of this relation will be called an {\it intersectional component}.
The {\it support} of a set $D$ of diagonals in a convex polygon is the set of vertices appearing in at least one of the diagonals of $D$. We denote this set by $Sup(D)$.
\end{defn}
\begin{exa}
The three intersectional components in Fig. \ref{Intersectional components} are $D_1=\big((1,7),(6,8)\big)$ (purple),$D_2=\big((1,5)\big)$ (green) and $D_3=\big((1,3),(2,4),(3,5)\big)$ (orange). The support of $D_1$ is the set $\{1,6,7,8\}$.
\end{exa}
We need the following:
\begin{lem}
\label{complete}
Let $G_n$ be a convex $n$-gon. If the graph of a diagonally framed dissection of $G_n$ contains a complete subgraph $K_k=(V_k, E_k)$ with $V_k=\{v_1,\dots,v_k\}$ and $(e,f)$ is a diagonal in the dissection that crosses one of the edges in $K_k$, then the graph with vertex set $V_k \cup \{e,f\}$ is a complete graph.
\end{lem}
\begin{proof}
We assume that $e,f \notin K_k$; the other cases can be treated by taking $V_k \setminus \{e,f\}$ instead.
Since $G_n$ is convex, the vertices $e$ and $f$ cannot be in the interior of $K_k$.  Since $(e,f)$ crosses at least one of the edges of $K_k$, the edge $(e,f)$ divides the vertices of the complete graph $K_k$ into two disjoint non-empty sets $A$ and $B$ such that $A\cup B=\{v_1,\dots,v_k\}$ and
each diagonal $(a,b)$ such that $a\in A$ and $b\in B$
crosses $(e,f)$ (see Fig. \ref{hexagon}). Since the dissection is diagonally framed,
the edges $(e,v_i),(v_i,f)$ exist in the dissection for each $1\leq i \leq k$, which completes the proof.\end{proof}

\begin{lem}\label{complete graph}
Let $G_n$ be the convex $n$-gon with a diagonally framed dissection.
Then for each intersectional component $C$ of the dissection of $G_n$, the graph induced by the support of $C$ is complete.
\end{lem}
\begin{proof}
We have to show that for each pair of vertices $a,b\in Sup(C)$, the edge $(a,b)$ is in $C$.
Let $a,b\in Sup(C)$. If $(a,b)\in C$, then we are done.
Otherwise, there are vertices $c,d \in Sup(C)$ such that $(a,c)\sim (b,d)$, i.e., there are vertices $\ell_1,\dots,\ell_k$ and $r_1,\dots,r_k$ such that $(\ell_i,r_i)$ crosses $(\ell_{i+1},r_{i+1})$ for $1 \leq i \leq k-1$ and such that $(\ell_1,r_1)=(a,c),(\ell_k,r_k)=(b,d)$.
We prove that the graph induced by the set of vertices $U_k=\{\ell_1,\dots,\ell_k,r_1,\dots,r_k\}$ is complete, and we are done.
We do this by induction on $k$: if $k=2$ then we have $U_2=\{\ell_1,\ell_2,r_1,r_2\}$ which induces a complete graph since $(\ell_1,r_1)$ crosses $(\ell_2,r_2)$ and the dissection of $G_n$ is diagonally framed by assumption.
Now let us assume that $U_k$ induces a complete graph and we have to prove that $U_{k+1}$ induces a complete graph. This is a direct consequence of Lemma \ref{complete}.\end{proof}
\begin{lem}
\label{2-polygon}

Let \(C=\{a=c_1<c_2<\cdots<c_m=b\}\) and \(D=\{a=d_1<d_2<\cdots<d_m=b\}\) be two distinct subpolygons of minimum order containing the edge \((a,b)\). Then some edge of \(C\) crosses some edge of \(D\)
\end{lem}

\begin{figure}
   \centering
  \includegraphics[scale=0.5]{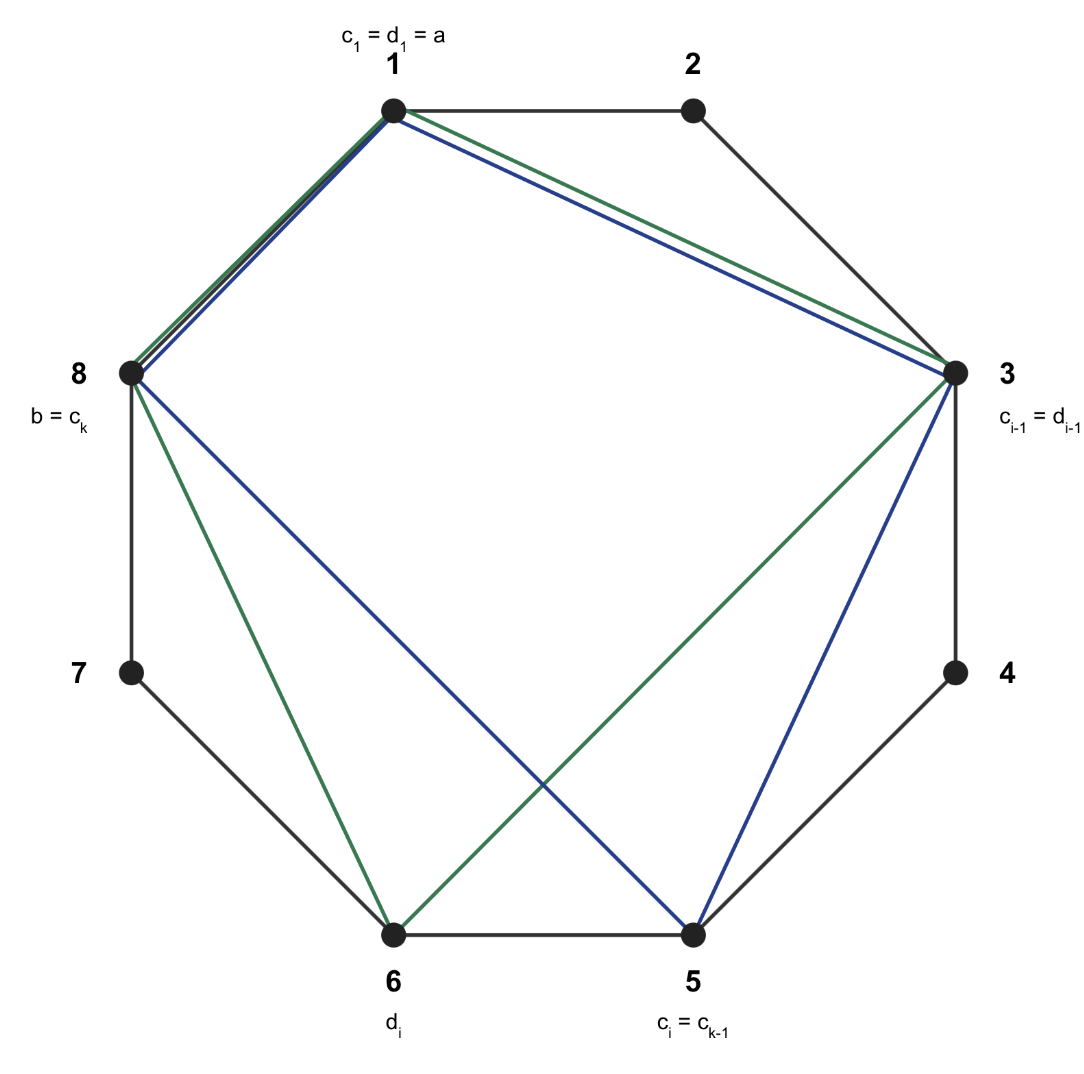}
    \caption{} 
   \label{distinct polygons}
\end{figure}

\begin{proof}
 Since $C\neq D$, there exists some minimal $i$ such that $c_i\neq d_i$. We may assume that $d_{i-1}<c_i<d_i$. Let $c_{k-1}$ be maximal among the vertices of $C$ such that $d_{i-1}<c_{k-1}<d_i$. It is clear that $c_k\neq d_i$, since otherwise the  boundary path
$c_{i-1},c_i,\ldots,c_k$
may be replaced by the shorter boundary path

$(c_{i-1}=d_{i-1
},c_k=d_i)$,
producing a subpolygon starting with $a$ and ending with $b$ with fewer
vertices than $C$.

Hence, $d_{i-1}<c_{k-1}<d_i<c_k$,
 so the edges
$(d_{i-1},d_i)$ and $(c_{k-1},c_k)$ cross. See Fig. \ref{distinct polygons} for an illustration. 
\end{proof}

Bouvel, Cioni, and Izart \cite{BCI} provided a formula for the number of interval posets with $n$ minimal elements.   They proved:
\begin{prop}
\label{threenon}
    The number of interval posets with $n>1$ minimal elements is
    $$\frac{1}{n}\sum\limits_{i=1}^{n-1}\sum\limits_{k=0}^{min\{i,\lfloor\frac{n-1-i}{2}\rfloor\}}\binom{n-1+i}{i} \binom{i}{k} \binom{n-2k-2}{i-1}.$$
\end{prop}
We are now ready to present the main result of this section.

\begin{figure}[!ht]
    \centering
   \includegraphics[scale=0.30]{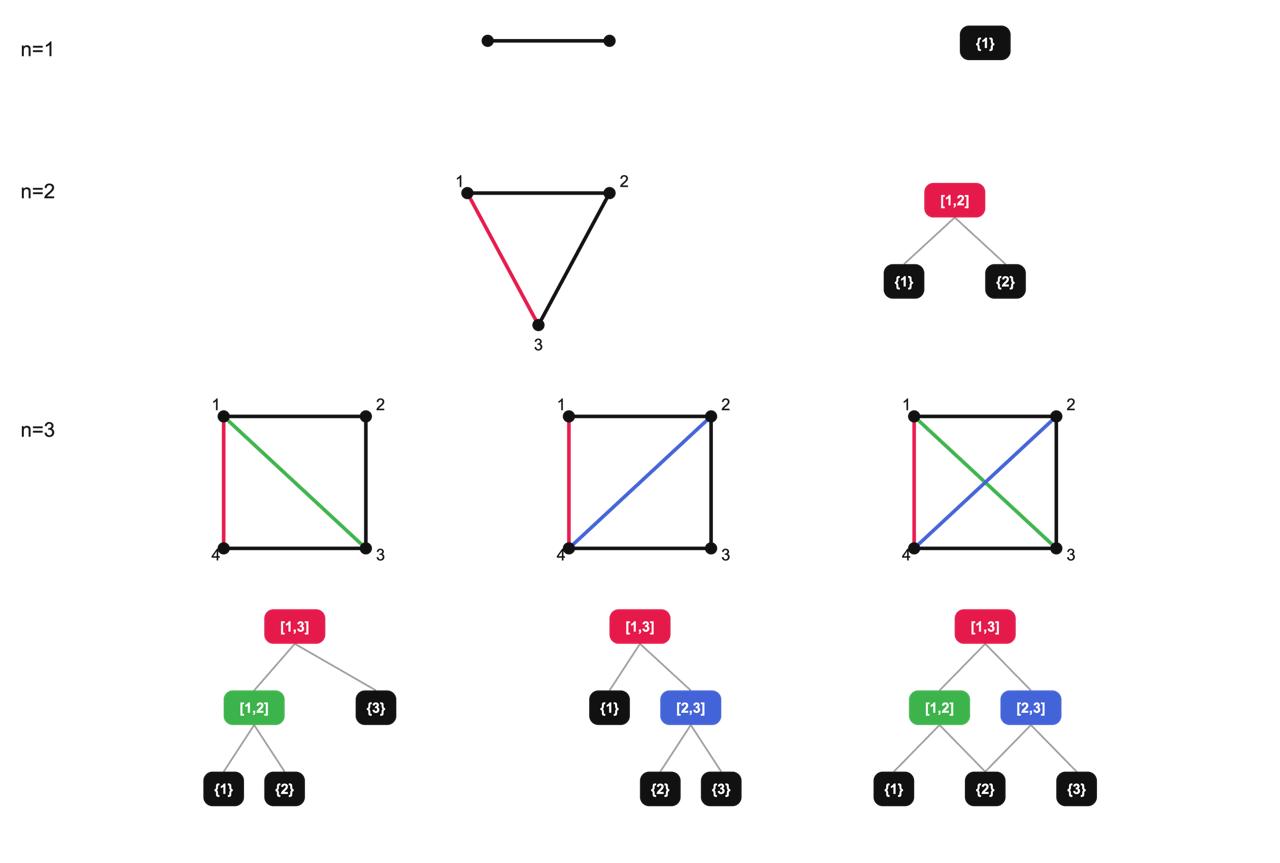}

    \caption{The bijection of Theorem \ref{bijection} for small values of $n$.}
    \label{all}
\end{figure}

\begin{thm}\label{bijection}
The mapping $\Phi$, from the set of interval posets with $n$ minimal elements to the set of diagonally framed dissections of the convex $(n+1)$-gon such that no quadrilateral is present, is bijective. Consequently, the sizes of these two sets are equal.
\label{bijection1}
\end{thm}
Fig. \ref{all} illustrates the bijection of Theorem \ref{bijection} for small values of $n$.

\begin{figure}[!ht]
    \centering   \includegraphics[scale=0.35]{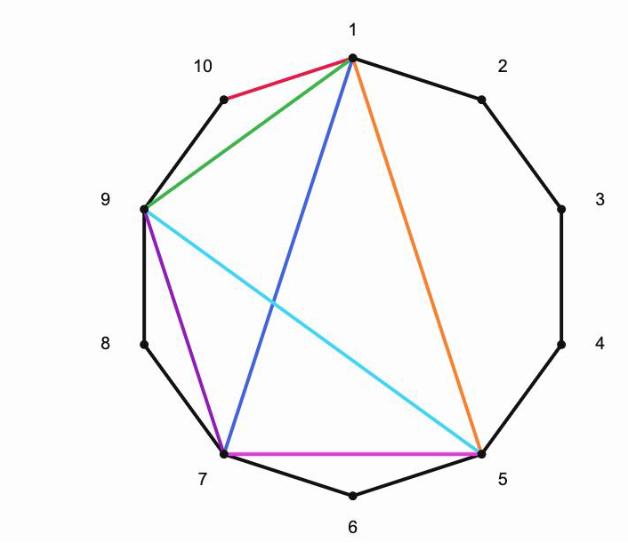}
    \caption{A dissection of a $10$-gon.}
    \label{example_argyle}
\end{figure}

\begin{figure}
 \includegraphics[scale=0.15]{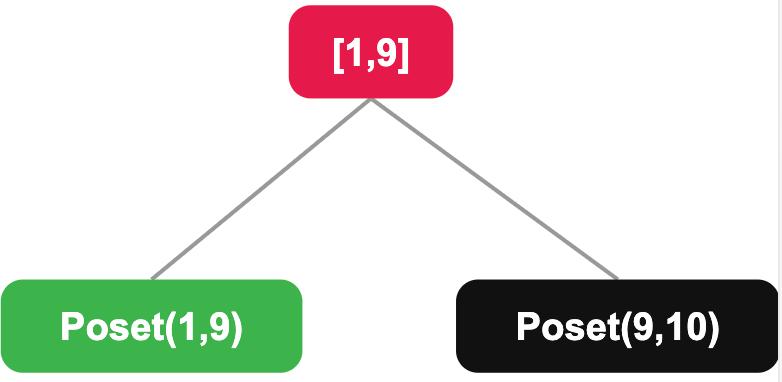}
 \caption{First step of Example \ref{from polygon to poset}}
 \label{step 1}
  \end{figure}

\begin{figure}
  \includegraphics[scale=0.15]{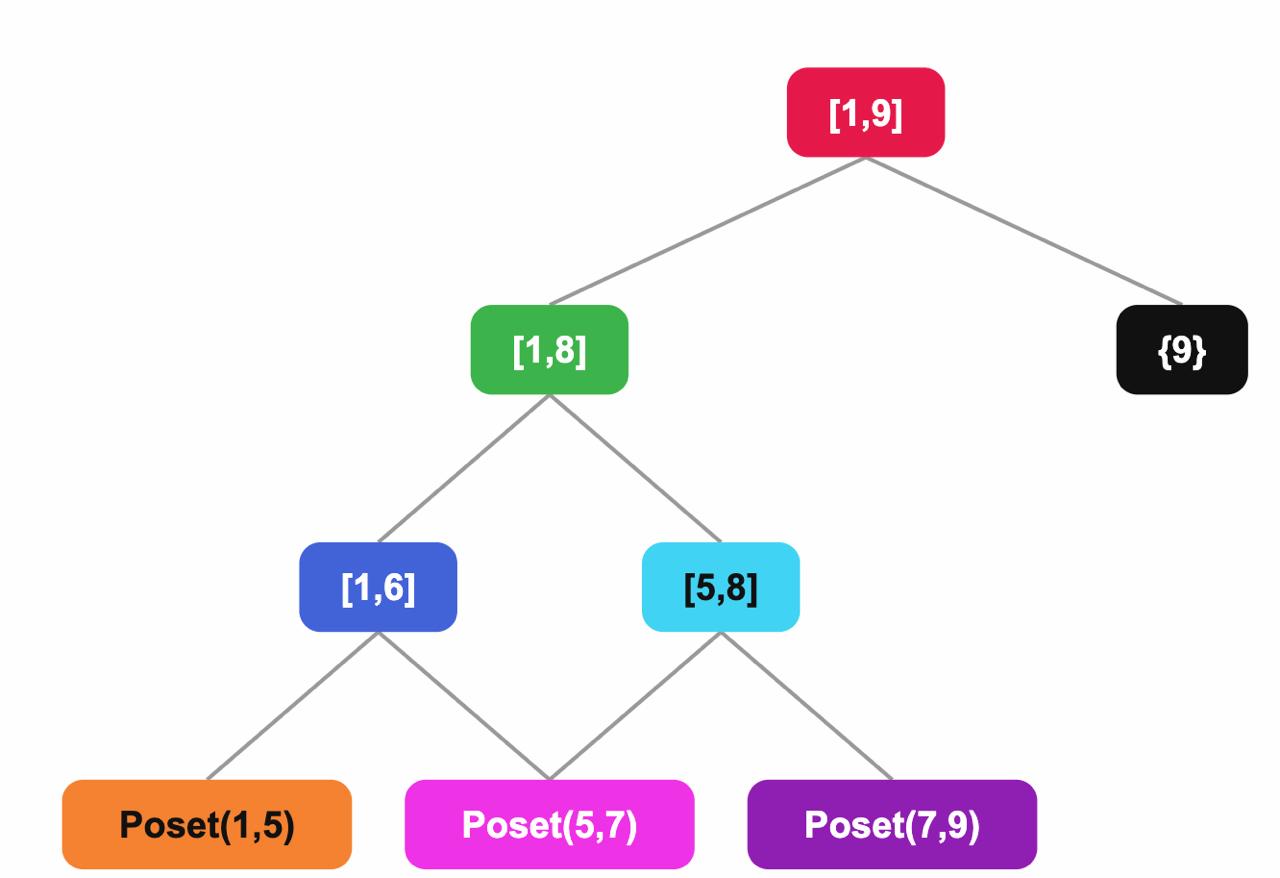}
  \caption{Second step of Example \ref{from polygon to poset}}
  \label{step 2}
  \end{figure}

\begin{figure}        \includegraphics[scale=0.15]{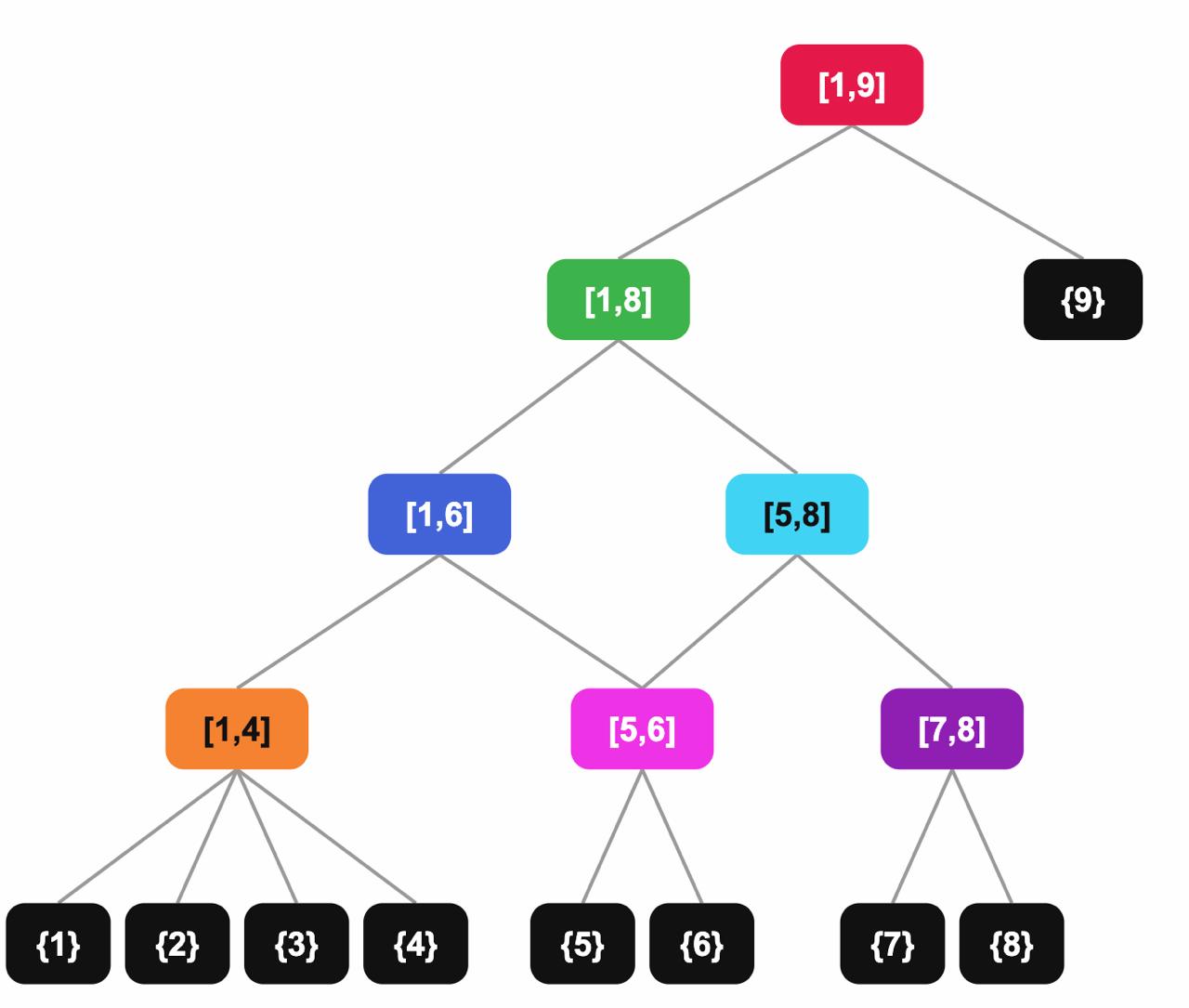}
    \caption{The third step of Example \ref{from polygon to poset}.}
    \label{step 3}
\end{figure}

\begin{proof}
The proof consists of two parts.
First, we show that every interval poset is mapped to a diagonally framed
dissection without quadrilaterals. 

Conversely, given such a dissection, we reconstruct the interval poset recursively using Theorem \ref{Interval}.

Let $P$ be an interval poset. By Lemma \ref{Phi(P) is diagonally framed}, $\Phi(P)$ is diagonally framed. Moreover, $\Phi(P)$ must not contain any quadrilateral, since  otherwise, if $a<b<c<d$ are such that $\{a,b,c,d\}$ is a quadrilateral then $P$ must contain the intervals $[a,b-1],[b,c-1],[c,d-1]$ and $[a,d-1]$ so the fourth interval covers exactly the first three. This is impossible by Proposition \ref{structure of interval posets 2}.

Conversely, given a diagonally framed dissection with no quadrilaterals of a convex $(n+1)$-gon with vertex set $V=\{1,\dots,n+1\}$, we match each edge $(a,b)$ of the dissection with the interval $[a,b-1]$, while if $a=b-1$ we get a minimal element $\{a\}$. It remains to show that the poset formed by these intervals is indeed an interval poset.

 We do this recursively, with the aid of Theorem \ref{Interval}.\\

We present a recursive algorithm for constructing the interval poset corresponding to a diagonally framed dissection. Given the surrounding edge $(1,n+1)$, the algorithm returns the entire interval poset with the $n$ minimal elements $
\{1\},
\dots,\{n\}.$

For each edge $(a,b)$ there are only two possibilities:

\begin{enumerate}[label=(\alph*)]
\item There is only one subpolygon of minimum order 
 $C=\{c_1,\dots ,c_\ell\}$, with 
$a=c_1<c_2<\cdots <c_\ell=b$ (see, for example, the edge $(1,5)$
in Fig. \ref{example_argyle}).

\item There are at least two different subpolygons of minimum order  $C$ as before and $D=\{d_1,\dots,d_{\ell}\}$  with $a=d_1<d_2<\cdots <d_\ell=b$, and then by Lemma \ref{2-polygon} and by Lemma \ref{complete graph}, the edge $(a,b)$ belongs to a
complete graph arising from an intersectional component  (see, for example, the edge $(1,9)$
in figure \ref{example_argyle}). 
\end{enumerate}

The algorithm is as follows:

$Poset(a,b)$\\
Consider the dissection of the subpolygon with vertex set $$V=\{a,a+1,\dots,b\}.$$
\begin{enumerate}
    \item If $b-a=1$ then
          return  the singleton $\{a\}$.
    \item Else, if we are in case (a) above, then return the poset having the interval $[a,b-1]$ as its maximal element and covering the $\ell-1$ elements $Poset(c_1,c_2),\dots, Poset(c_{\ell-1},c_\ell)$. According to \ref{argyle}, this is either an argyle poset of order 2 or a dual claw poset.
    \item Else, if we are in case (b) then return 
    an argyle poset with $\ell-1$ minimal elements, having the interval $[a,b-1]$ as its maximal element and containing the intervals $[c_i,c_j-1]$ for each $1\leq i<j-1\leq \ell$.
    Finally, for each minimal element of the argyle, apply
    $Poset(c_i,c_{i+1})$ for $1\leq i\leq \ell-1$.
\end{enumerate}

By Theorem \ref{Interval}, we are done. \end{proof}
\begin{exa}
\label{from polygon to poset}
Consider the dissection of the $10$-gon depicted in Fig.~\ref{example_argyle}.
We begin by applying $Poset(1,10)$. 

Since the unique subpolygon of minimum order 
containing the edge $(1,10)$ is the triangle
$C=\{1,9,10\}$, step~(2) of the algorithm yields the poset shown in Fig.~\ref{step 1}.

Next, we apply $Poset(9,10)$, which returns the singleton $\{9\}$ by
step~(1). 

We also apply $Poset(1,9)$. Note that there are two different subpolygons of minimum order containing the edge $(1,9)$: $C=\{1,5,9\}$ and $D=\{1,7,9\}$. Hence, by step~(3), we obtain the argyle poset corresponding to the complete graph on the
vertex set $\{1,5,7,9\}$ which is presented in 
Fig.~\ref{step 2}.

We now apply $Poset(1,5)$. According to step~(2), this yields a dual claw
poset whose maximal element is the interval $[1,4]$ and which covers
$Poset(1,2)$, $Poset(2,3)$, $Poset(3,4)$, and $Poset(4,5)$. These, in
turn, return the singletons $\{1\}$, $\{2\}$, $\{3\}$, and $\{4\}$,
respectively.

Similarly, $Poset(5,7)$ yields the interval $[5,6]$, which covers the
singletons $\{5\}$ and $\{6\}$, while $Poset(7,9)$ yields the interval
$[7,8]$, which covers the singletons $\{7\}$ and $\{8\}$. The final
interval poset is shown in Fig.~\ref{step 3}.
\end{exa}


\section{Tree interval posets}
Recall that a tree poset is a poset whose Hasse diagram is a tree.
\begin{figure}
    \centering
   \includegraphics[scale=0.35]{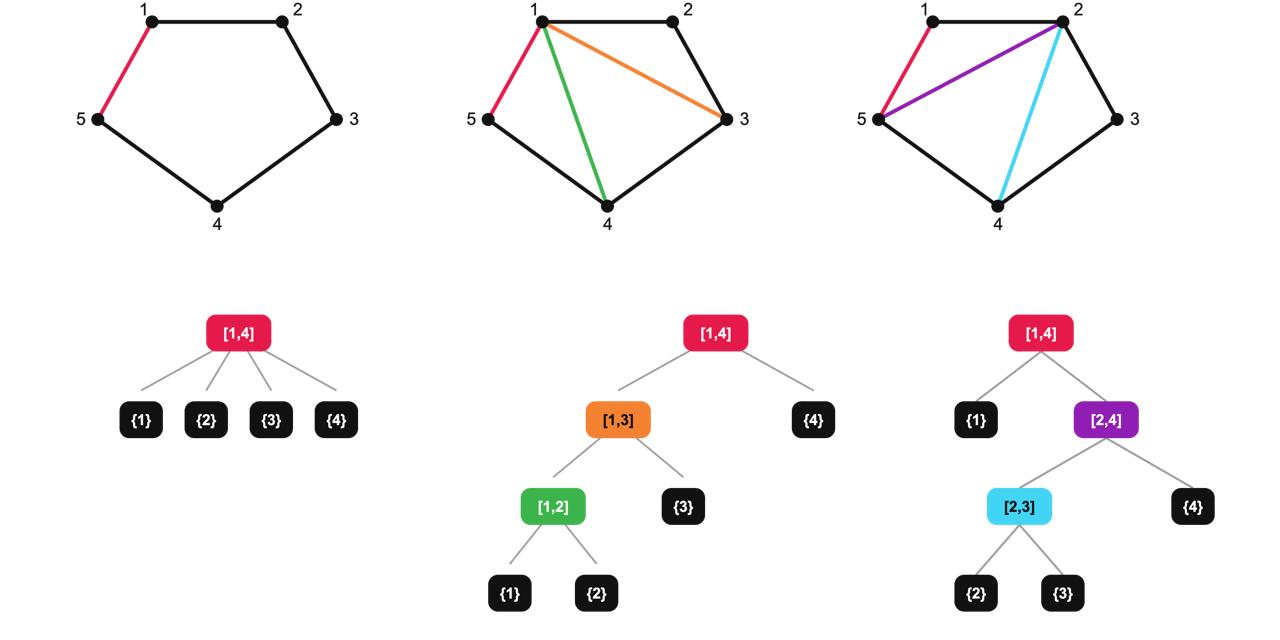}
    \includegraphics[scale=0.35]{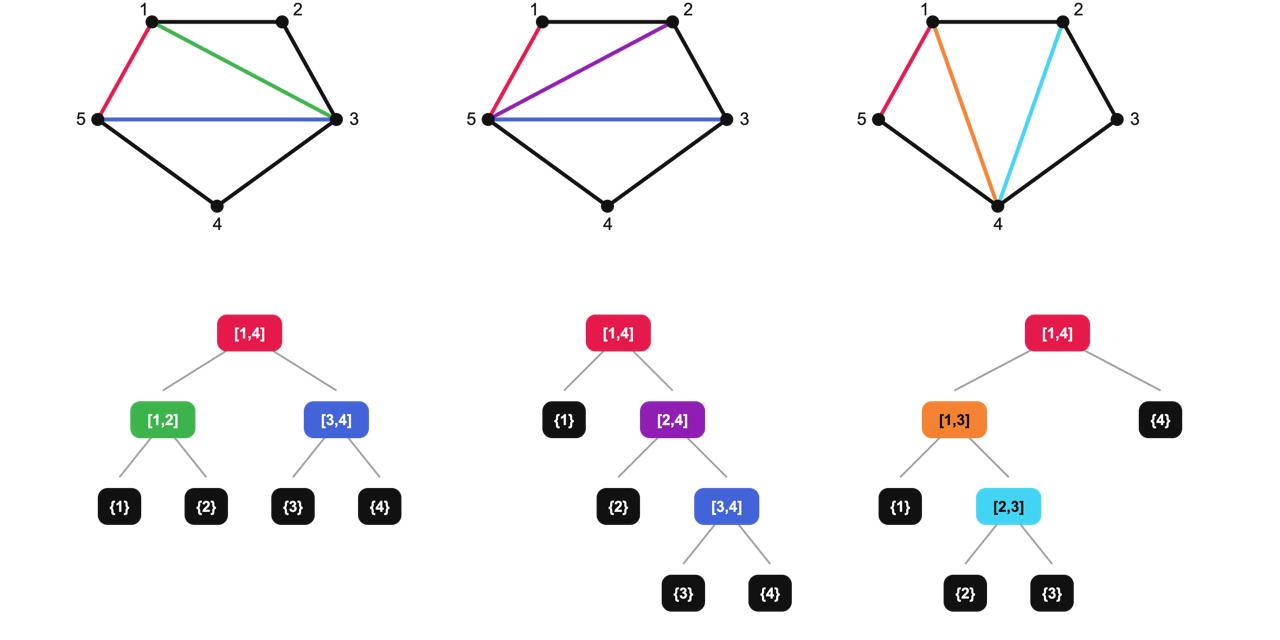}
   \caption{Examples of the bijection for tree interval posets for $n=4$.}
    \label{trees}
\end{figure}
In \cite{BCI}, the authors calculated the number of tree interval posets using generating functions and mentioned that this is equal to the number of ways to place non-crossing diagonals in a convex $(n+1)$-gon such that no quadrilateral is created (sequence A054515 from OEIS \cite{OEIS}).
Using the function $\Phi$ defined above, we obtain a combinatorial proof of this result.
\begin{thm}\label{enumerate trees}
The restriction of the mapping $\Phi$, defined above, to the set of tree interval posets with $n$ minimal elements is a bijection to the set of non-crossing dissections of the convex $(n+1)$-gon such that no quadrilateral is present. Consequently, the sizes of these two sets are equal.
\end{thm}
Figure \ref{trees} illustrates the bijection of Theorem \ref{enumerate trees} for $n=4$.

\begin{proof}

By Theorem \ref{bijection}, the mapping $\Phi$ is a bijection between interval
posets with $n$ minimal elements and diagonally framed dissections of a
convex $(n+1)$-gon containing no quadrilateral.

Note that the intervals $[a,b]$ and $[c,d]$ intersect and are incomparable if and only if the corresponding diagonals $(a,b+1)$ and $(c,d+1)$ cross each other (i.e., intersect in an interior point).  

 It therefore remains to
show that an interval poset $P$ is a tree if and only if it does not have any incomparable intersecting intervals. 

Suppose first that $P$ has two incomparable intersecting intervals $I$ and $J$.  
 By Observation \ref{structure of interval posets 1}, both
$I\cap J$ and $I\cup J$ are intervals in $P$.
Hence there
is a saturated chain from $I\cap J$ to $I\cup J$ passing through $I$,
and another saturated chain passing through $J$. Since $I$ and $J$ are
incomparable, these two chains must diverge and subsequently rejoin.
Their union therefore contains a cycle in the Hasse diagram of $P$.
Consequently, $P$ is not a tree. See Fig. \ref{general diamond} with $I=[1,4],J=[3,6], I\cap J=[3,4]$, and $I\cup J=[1,6]$.

Conversely, suppose that no two
incomparable intervals of $P$ intersect. In particular, no argyle poset
of order greater than $2$ occurs in the recursive decomposition of $P$ 
given by Theorem \ref{Interval}; otherwise, one would have two intervals covering a common interval.  
Thus $P$ is constructed recursively using only
dual claw posets and argyle posets of order $1$ or $2$, whose Hasse diagrams
are trees. Since replacing a minimal element of a tree by either of
these posets preserves the tree property, the Hasse diagram of $P$ is
a tree.

Therefore, $P$ is a tree if and only if $\Phi(P)$ is non-crossing, and
the restriction of $\Phi$ gives the required bijection.

\end{proof}

\begin{figure}

\begin{tikzpicture}[
  x=1.55cm,y=1.35cm,
  interval/.style={draw=black!75, fill=white, rounded corners=1.2pt,
                   inner xsep=4pt, inner ysep=2.5pt, font=\small},
  diamond interval/.style={interval, draw=red!75!black, fill=red!3,
                           very thick},
  cover/.style={draw=black!45, line width=.45pt},
  diamond cover/.style={draw=red!75!black, line width=1.15pt}
]
  \node[interval] (i11) at (0,0) {$[1]$};
  \node[interval] (i22) at (1,0) {$[2]$};
  \node[interval] (i33) at (2,0) {$[3]$};
  \node[interval] (i44) at (3,0) {$[4]$};
  \node[interval] (i55) at (4,0) {$[5]$};
  \node[interval] (i66) at (5,0) {$[6]$};
  \node[interval] (i12) at (0.5,1) {$[1,2]$};
  \node[interval] (i23) at (1.5,1) {$[2,3]$};
  \node[diamond interval] (i34) at (2.5,1) {$[3,4]$};
  \node[interval] (i45) at (3.5,1) {$[4,5]$};
  \node[interval] (i56) at (4.5,1) {$[5,6]$};
  \node[interval] (i13) at (1,2) {$[1,3]$};
  \node[diamond interval] (i24) at (2,2) {$[2,4]$};
  \node[diamond interval] (i35) at (3,2) {$[3,5]$};
  \node[interval] (i46) at (4,2) {$[4,6]$};
  \node[diamond interval] (i14) at (1.5,3) {$[1,4]$};
  \node[interval] (i25) at (2.5,3) {$[2,5]$};
  \node[diamond interval] (i36) at (3.5,3) {$[3,6]$};
  \node[diamond interval] (i15) at (2,4) {$[1,5]$};
  \node[diamond interval] (i26) at (3,4) {$[2,6]$};
  \node[diamond interval] (i16) at (2.5,5) {$[1,6]$};

  \draw[cover] (i12) -- (i11);
  \draw[cover] (i12) -- (i22);
  \draw[cover] (i23) -- (i22);
  \draw[cover] (i23) -- (i33);
  \draw[cover] (i34) -- (i33);
  \draw[cover] (i34) -- (i44);
  \draw[cover] (i45) -- (i44);
  \draw[cover] (i45) -- (i55);
  \draw[cover] (i56) -- (i55);
  \draw[cover] (i56) -- (i66);
  \draw[cover] (i13) -- (i12);
  \draw[cover] (i13) -- (i23);
  \draw[cover] (i24) -- (i23);
  \draw[cover] (i35) -- (i45);
  \draw[cover] (i46) -- (i45);
  \draw[cover] (i46) -- (i56);
  \draw[cover] (i14) -- (i13);
  \draw[cover] (i25) -- (i24);
  \draw[cover] (i25) -- (i35);
  \draw[cover] (i36) -- (i46);
  \draw[cover] (i15) -- (i25);
  \draw[cover] (i26) -- (i25);
  \draw[diamond cover] (i24) -- (i34);
  \draw[diamond cover] (i35) -- (i34);
  \draw[diamond cover] (i14) -- (i24);
  \draw[diamond cover] (i36) -- (i35);
  \draw[diamond cover] (i15) -- (i14);
  \draw[diamond cover] (i26) -- (i36);
  \draw[diamond cover] (i16) -- (i15);
  \draw[diamond cover] (i16) -- (i26);
\end{tikzpicture}
\caption{An example of a diamond subposet inside the interval poset of $123456$}
\label{general diamond}
\end{figure}

\section{Interval posets of block-wise simple permutations}
In \cite{BERS}, the authors of this paper introduced the notion of block-wise simple permutations. We cite here the definition:
\begin{defn}\label{def by 2 blocks}
A permutation $\pi \in \mathcal{S}_n$ is called {\em block-wise simple} if it has no interval of the form $p_1\oplus p_2$ or $p_1 \ominus p_2$,  where $\oplus$ and $\ominus$ stand for direct and skew sums of permutations, respectively. See \cite{BE} for definitions of direct and skew sums. 
\end{defn}
There are no block-wise simple permutations of sizes $2$ and $3$.
For $n \in \{4,5,6\}$, a permutation is block-wise simple if and only if it is simple. The first nontrivial examples of block-wise simple permutations are of size 7, for instance $4253716$.

In \cite{BERS}, the authors enumerated the interval posets of block-wise simple permutations, using methods 
 similar to the ones appeared in \cite{BCI}, to prove the following proposition, which appears in \cite{BERS} as Theorem 4.4. 
\begin{prop}
The number of interval posets corresponding to block-wise simple permutations of size $n\geq 4$ is
$$\frac{1}{n}\sum_{i=1}^{\lfloor\frac{n-1}{3}\rfloor}{n+i-1 \choose i}{n-2i-2 \choose i-1}.$$
\end{prop}
 This is sequence A054514 from OEIS \cite{OEIS}, which also counts the number of ways to place non-crossing diagonals in a convex $(n+1)$-gon such that there are no triangles or quadrilaterals. We now provide a combinatorial proof of this claim, using the bijection defined above.
\begin{figure}
    \centering
     \includegraphics[scale=0.20]{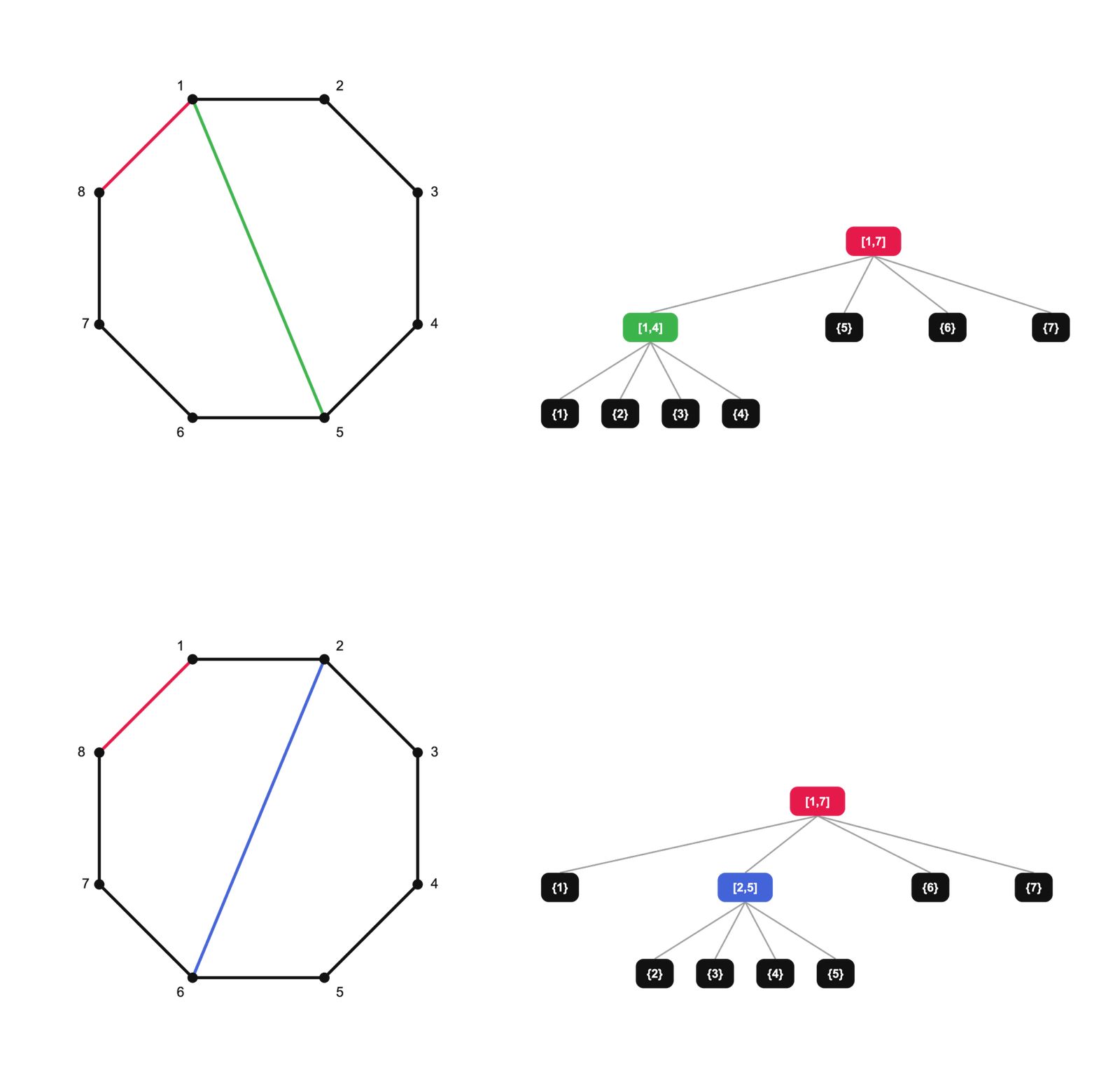}
    \includegraphics[scale=0.20]{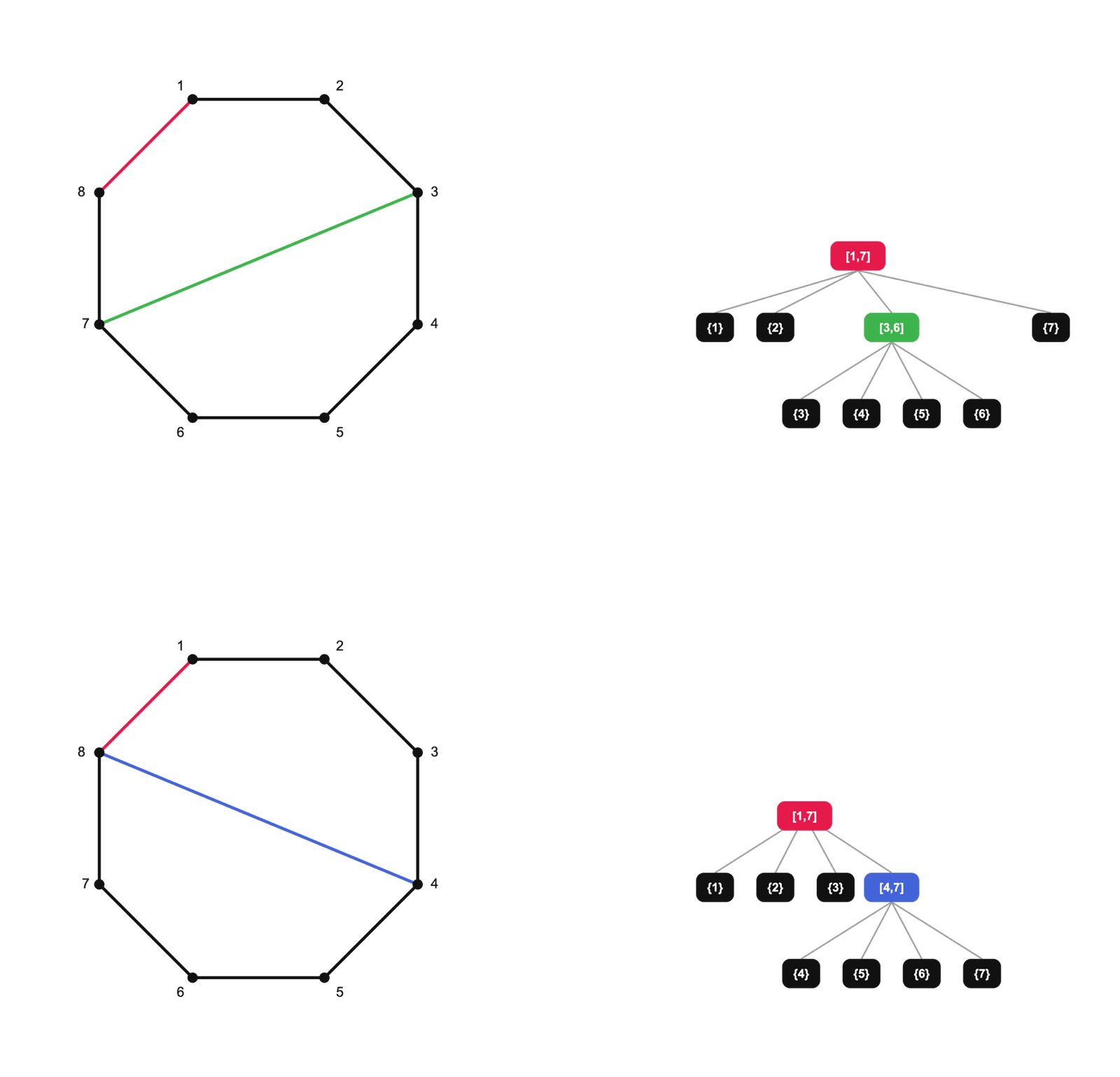}
   \caption{Examples of the bijection for block-wise simple interval posets. }
    \label{blockwise}
\end{figure}
\begin{thm}\label{bijection block-wise}
The restriction of the mapping $\Phi$  to the set of interval posets of block-wise simple permutations of size $n$ is a bijection to the set of non-crossing dissections of a convex $(n+1)$-gon such that no triangles or quadrilaterals are present. Consequently, the sizes of these sets are equal.
\end{thm}
Fig. \ref{blockwise} illustrates the bijection of Theorem \ref{bijection block-wise} for $n=7$.
\begin{proof}
In \cite[Theorem ~6.3]{T}, the author proved that $P=P(\pi)$ is a tree interval poset if and only if $\pi$ contains no interval of the form $p_1\oplus p_2\oplus p_3$ or $p_1\ominus p_2 \ominus p_3$.
From this, and by Definition \ref{def by 2 blocks}, it follows that an interval poset of a block-wise simple permutation is a tree. Hence, by Theorem \ref{enumerate trees}, $\Phi(P)$ has no crossing diagonals and has no quadrilaterals. Suppose, to the contrary, that $\Phi(P)$ contains a triangle with edges $(a,b)$,
$(b,c)$, and
$(a,c)$ where
$a<b<c$. Then $P$ contains the intervals $[a,b-1]$, $[b,c-1]$, and $[a,c-1]$. In this case $[a,c-1]$ covers both $[a,b-1]$ and $[b,c-1]$
which means that $\pi$ contains a sum $p\oplus q$ or a skew sum $p \ominus q$, where $p$ is the part of $\pi$ containing the elements of the interval $[a,b-1]$ and $q$ is the part of $\pi$ containing the elements of the interval $[b,c-1]$, contradicting the definition of block-wise simple permutations.

Conversely, assume that $\Phi(P)$ has no quadrilaterals, no crossing diagonals and no triangles.

Thus, by Theorem \ref{enumerate trees}, \(P\) is a tree, and its recursive decomposition contains only dual claw posets and argyle posets of order two.

Since $\Phi(P)$ has no triangles, $P$ has no argyle of order two. (Indeed, if $P$ has an argyle of order two, then  $P$ contains an interval $[a,c]$ which covers the intervals 
$[a,b],[b+1,c]$ for some $a<b<c$ and therefore $\phi(P)$ contains the diagonals $(a,b+1),(b+1,c+1)$ and $(a,c+1)$ which produce a triangle). As a result, $P$ contains only dual claw subposets and thus, according to definition \ref{dual claw argyle},
every non-minimal element of $P$ covers more than $2$ elements. Therefore $\pi$ has no intervals of the form $p_1\oplus p_2$ or $p_1 \ominus p_2$ so $P$ is an interval poset of a block-wise simple permutation.   
\end{proof}

\section{Interval posets of separable permutations}\label{Interval posets of separable permutations}
A permutation is called {\it separable} if it can be obtained from the trivial permutation $1$ by direct and skew sums. More explicitly, $1$ is separable and a permutation $\pi\neq 1$ is separable if it can be written as $\pi=\pi_1 \oplus \pi_2$ or $\pi=\pi_1\ominus \pi_2$ such that $\pi_1$ and $\pi_2$ are themselves separable.
 For example, $\pi=45213867=\pi_1 \oplus \pi_2$ where $\pi_1=45213$ and $\pi_2=312$. Now, $\pi_1=12\ominus 213$ and $\pi_2=1\ominus (1\oplus 1)$. Altogether, we obtain $$\pi=((1\oplus 1)\ominus ((1\ominus 1)\oplus 1))\oplus(1 \ominus(1\oplus 1)).$$
The notion of separable permutations is well studied (see, for example, \cite{K}, \cite{AN}, and \cite{BBL}, where they are characterized by avoiding the forbidden patterns $2413$ and $3142$).
Recall that an interval poset is called {\em binary} if no element covers more than two elements. Note that a binary interval poset is not necessarily a tree.  By Theorem 6.2 in \cite{T},
a permutation $\pi$ is separable if and only if its interval poset is binary.
The number of separable permutations of size $n$ is known to be the (large) Schröder number (sequence A006318 in the OEIS \cite{OEIS}).
In \cite[Theorem ~6.3]{T}, the author proved that each interval poset of a separable permutation corresponds to exactly two separable permutations, from which she concluded that the number of such interval posets is half this number, namely the small Schröder number (sequence A001003 in the OEIS \cite{OEIS}).
It is known (\cite{S}) that small Schröder numbers also count the dissections of the $(n+1)$-gon such that crossing diagonals are prohibited. Here, we present a bijection, different from the one we used earlier in this paper, to prove this connection in a combinatorial way.
\begin{figure}
    \centering
     \includegraphics[scale=0.25]{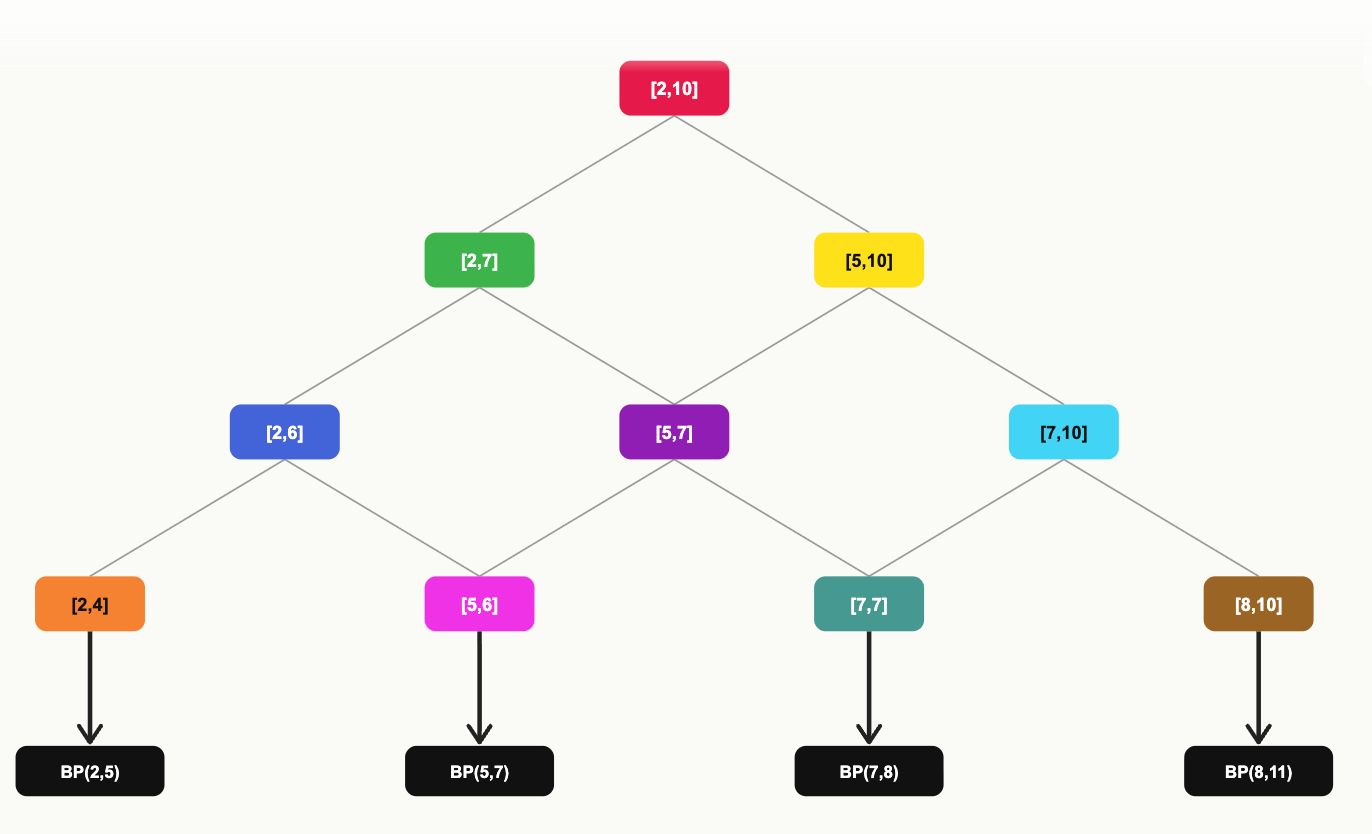}
    \includegraphics[scale=0.15]{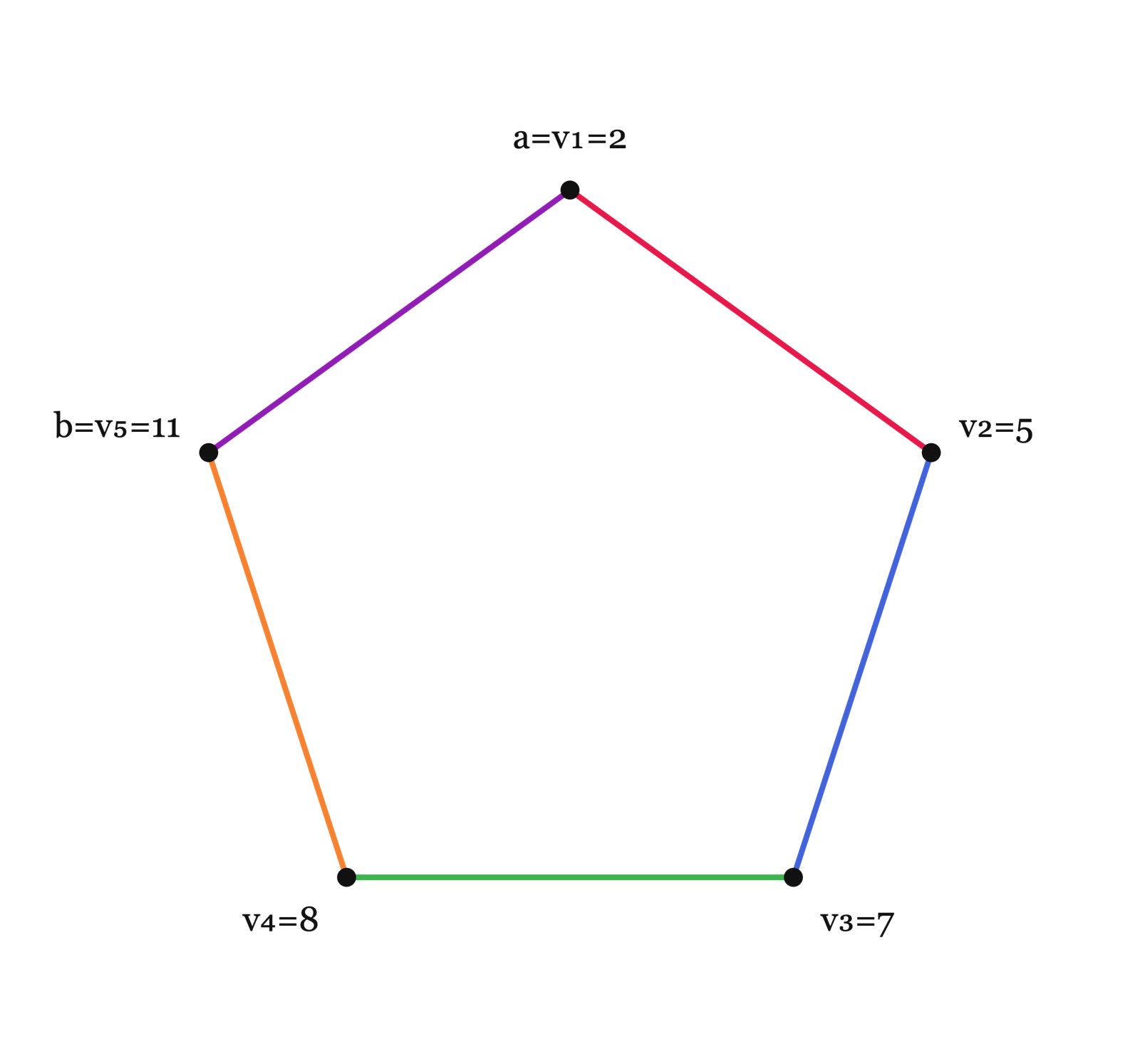}
   \caption{Close-up of case (b) of the algorithm BinPoset}
    \label{zoom up}
\end{figure}

\begin{figure}[!h]
    \centering
   \includegraphics[scale=0.25]{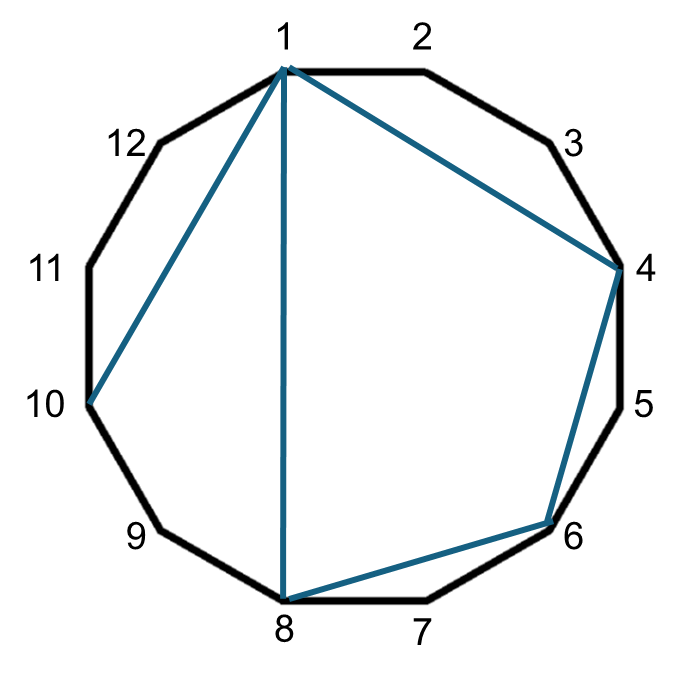}
   \caption{Example of the bijection for the interval poset of separable permutations.}
    \label{example for bijection separable}
\end{figure}

\begin{figure}[!h]
   \includegraphics[scale=0.35]{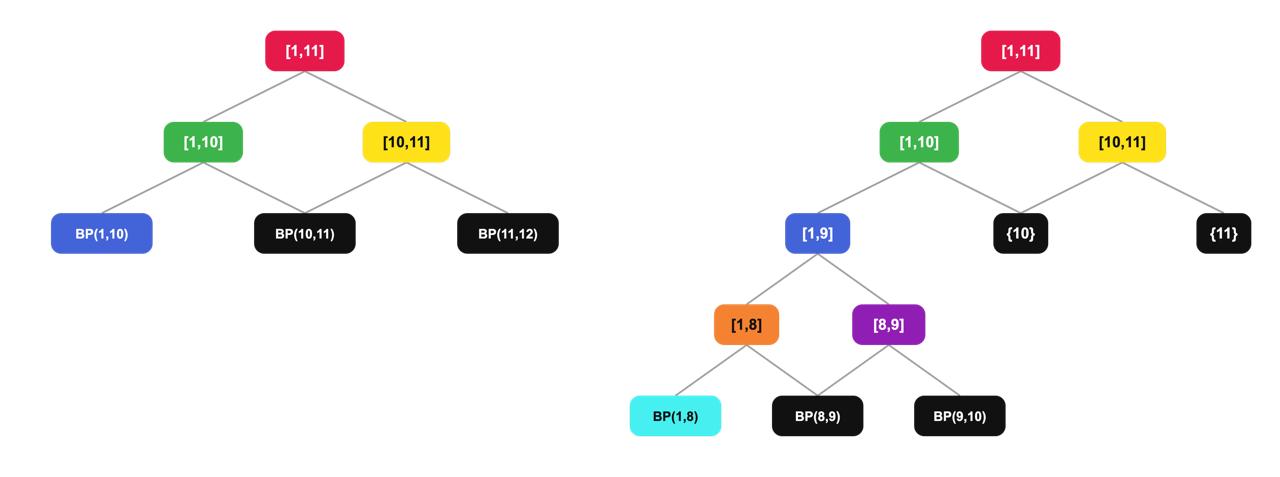}
   \caption{From left to right: the first and the second step in Example \ref{example separable}.}
    \label{BP2}
\end{figure}

\begin{figure}[!h]
   \includegraphics[scale=0.65]{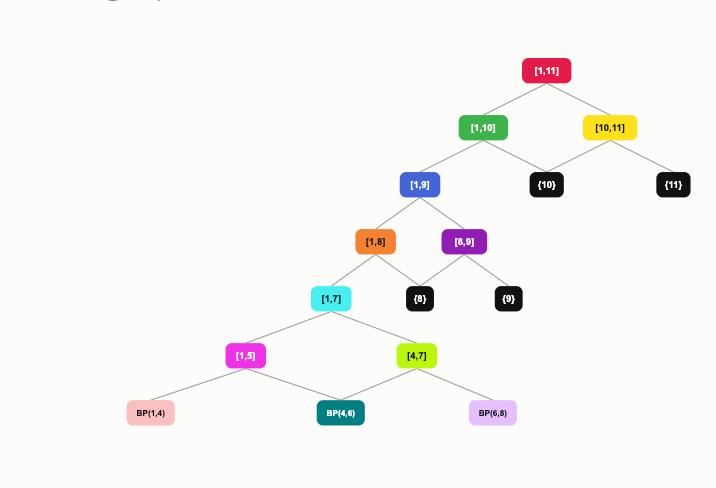}
   \caption{ The third step in Example \ref{example separable}.}
    \label{BP3}
\end{figure}

\begin{figure}[!h]
   \includegraphics[scale=0.65]{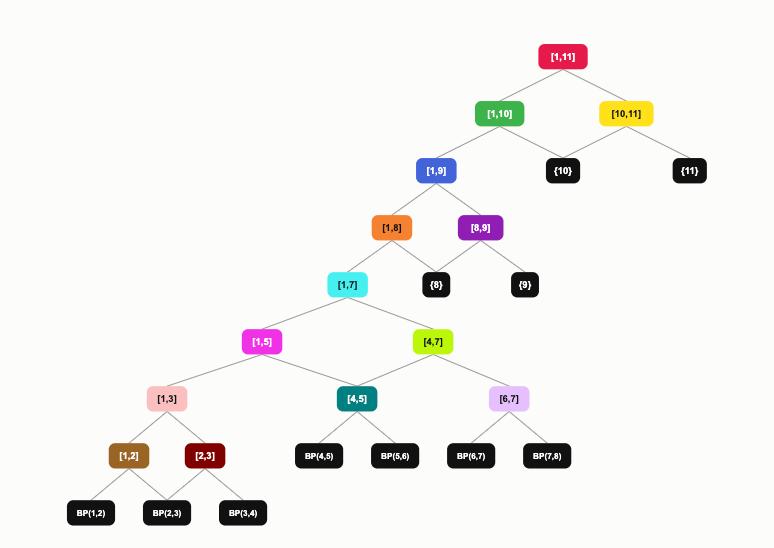}
   \caption{ The fourth step in Example \ref{example separable}.}
    \label{BP4}
\end{figure}
\begin{figure}[!h]
    \centering
   \includegraphics[scale=0.35]{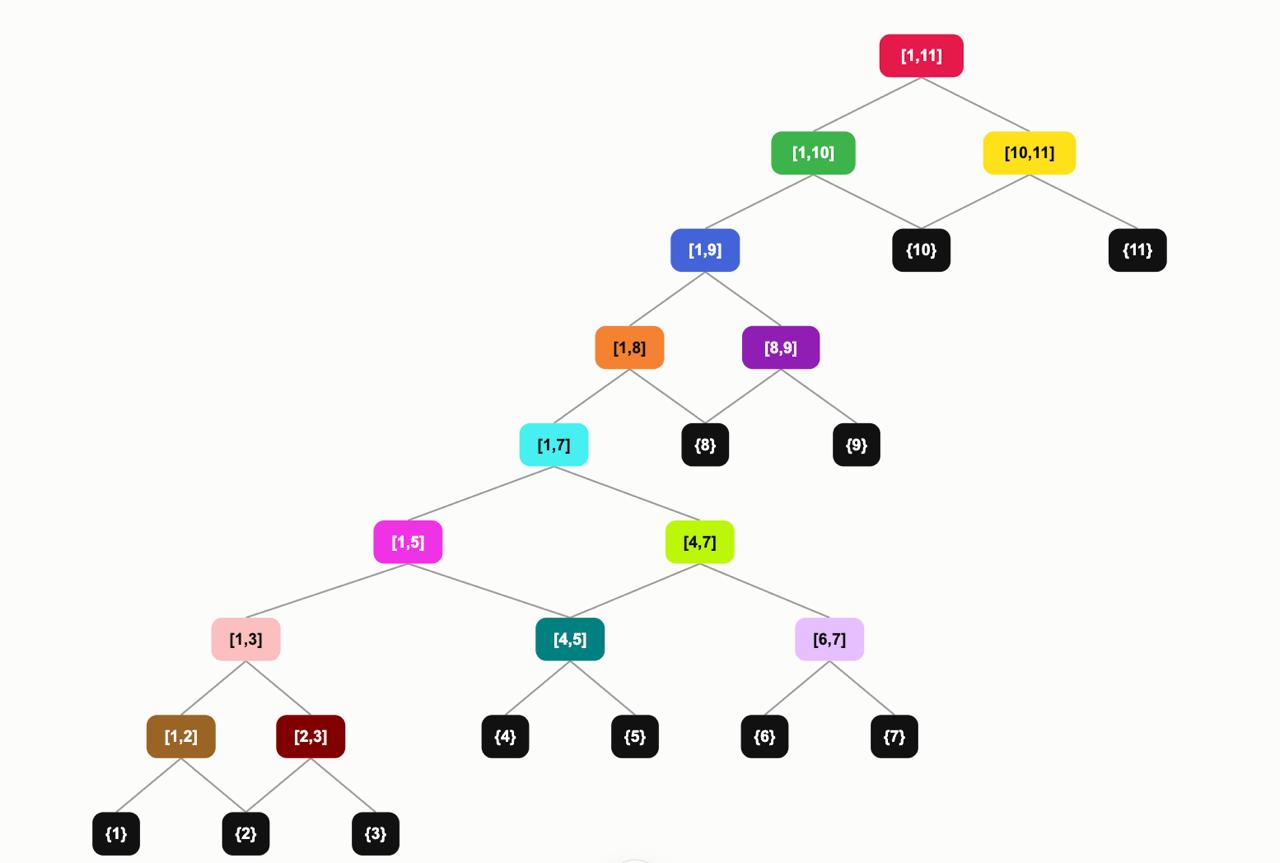}
   \caption{The final step in Example \ref{example separable}.}
    \label{BP5}
\end{figure}

\begin{thm}
The set of interval posets corresponding to separable permutations of size $n$ is in bijection with the set of dissections of the $(n+1)$-gon such that crossing diagonals are prohibited.  Consequently, the sizes of these sets are equal.
\end{thm}
\begin{proof}
We establish a mapping from the set of dissections of the $(n+1)$-gon, having no crossing diagonals, to the set of interval posets of separable permutations, i.e., the binary interval posets.
We present a recursive algorithm $BinPoset(a,b)$ which, given the edge $(1,n+1)$, returns the corresponding interval poset with $n$ minimal elements. \\

$BinPoset(a,b)$\\
Consider the dissection of the subpolygon with vertex set $$V=\{a,a+1,\dots,b\}.$$
\begin{enumerate}
    \item If $b-a=1$, then
          return the singleton $\{a\}$.
    \item Else, there is a subpolygon of minimal order $k$ that contains the edge $(a,b)$. We now have two options:
\begin{enumerate}
\item If $k=3$, then denote the set of vertices of the triangle by $\{a,m,b\}$ with $a<m<b$. Return a poset, with $[a,b-1]$ as its maximal element, which covers the elements $BinPoset(a,m)$ and $BinPoset(m,b)$.
\item Otherwise, if $k>3$, then let the vertices of the subpolygon be $\{a=v_1,\dots, v_k=b\}$ with $v_1<v_2\cdots< v_k$.
Return an argyle poset with  $[a,b-1]$ as its root, and for each $1 \leq i<j-1 \leq k$ include the interval $[v_i,v_j-1]$ with its two descendants: $[v_i,v_{j-1}-1]$ and $[v_{i+1},v_j-1]$. Instead of the minimal elements, we include the elements $BinPoset(v_i,v_{i+1})$ for $1 \leq i \leq k-1$. Fig. \ref{zoom up} is a close-up of this case.

\end{enumerate}
    \end{enumerate}

On the other hand, every binary interval poset is a composition of argyle posets. Given such an interval poset with \(n\) minimal elements, construct a dissection of the convex \((n+1)\)-gon as follows: for each argyle in the interval poset, draw the diagonal \((a,b+1)\), where \([a,b]\) is its maximal element, except when \([a,b]=[1,n]\), the maximal element of the entire poset.
The interval \([1,n]\) is excluded because it corresponds to the outer edge \((1,n+1)\) of the polygon.

\end{proof}

\begin{exa}
\label{example separable}
    Given the dissection of the $12$-gon depicted in Fig. \ref{example for bijection separable}, we construct the binary interval poset step by step. 
    
    First, we consider the edge $(1,12)$ which is contained in the subpolygon of minimum order $\{1,10,11,12\}$. According to part $(2) (b)$ of the algorithm,  we get the argyle poset that can be seen on the left-hand side of Fig. \ref{BP2}.
    
Next, according to case $(1)$, $BinPoset(11,12)$ is replaced by the singleton $\{11\}$, $BinPoset(10,11)$ is replaced by the singleton $\{10\}$, and according to $(2)(b)$, the call $BinPoset(1,10)$ produces an argyle poset with $[1,9]$ as its maximal element.  The result can be seen on the right side of Fig. \ref{BP2}.

In the next step, $BinPoset(8,9)$ and $BinPoset(9,10)$ are replaced by the singletons $\{8\}$ and $\{9\}$ respectively, and since $(1,8)$ is contained in the subpolygon $\{1,4,6,8\}$, the call $BinPoset(1,8)$ produces an argyle with maximal element $[1,7]$ as in Fig. \ref{BP3}.
The two remaining steps are depicted in Fig. \ref{BP4} and in Fig. \ref{BP5}.

\end{exa}
\section*{Acknowledgments}
We would like to thank the anonymous referees for their thorough review of the first version of this paper and their helpful suggestions.

\end{document}